\documentclass[11pt]{article}
\usepackage{amsmath,amssymb,amsthm}
\usepackage{dsfont}
\usepackage{graphicx}
\usepackage{geometry}
\usepackage[british]{babel}
\usepackage[utf8]{inputenc}
\usepackage{xcolor}
\usepackage[linesnumbered,ruled,vlined,commentsnumbered]{algorithm2e}
\usepackage{cite}
\usepackage{url}
\usepackage{diagbox}
\usepackage{hyperref}
\usepackage{authblk}
\usepackage{booktabs}
\usepackage{todonotes}
\usepackage{caption}
\usepackage{comment}
\def\epsilon{\varepsilon}

\DeclareCaptionLabelFormat{AppendixTables}{A.#2}

\DeclareMathOperator{\argmin}{argmin}

\newtheorem{definition}{Definition}

\allowdisplaybreaks

\begin{document}
\title{Asymptotic-Preserving Neural Networks for multiscale hyperbolic models of epidemic spread}

\date{}

\author[1,2]{Giulia Bertaglia}
\author[3]{Chuan Lu}
\author[2]{Lorenzo Pareschi}
\author[3]{Xueyu Zhu}
\affil[1]{Istituto Nazionale di Alta Matematica “Francesco Severi”,\newline P.le Aldo Moro 5, 00185 Roma, Italy.}
\affil[2]{Department of Mathematics and Computer Science, University of Ferrara,\newline Via Machiavelli 30, 44121 Ferrara, Italy.}
\affil[3]{Department of Mathematics, University of Iowa,\newline Iowa City, IA 52246, USA.}
\maketitle
\begin{abstract}
When investigating epidemic dynamics through differential models, the parameters needed to understand the phenomenon and to simulate forecast scenarios require a delicate calibration phase, often made even more challenging by the scarcity and uncertainty of the observed data reported by official sources. 
In this context, Physics-Informed Neural Networks (PINNs), by embedding the knowledge of the differential model that governs the physical phenomenon in the learning process, can effectively address the inverse and forward problem of data-driven learning and solving the corresponding epidemic problem.
In many circumstances, however, the spatial propagation of an infectious disease is characterized by movements of individuals at different scales governed by multiscale PDEs. This reflects the  heterogeneity  of a region or territory in relation to the dynamics within cities and in neighboring zones.
In presence of multiple scales, a direct application of PINNs generally leads to poor results due to the multiscale nature of the differential model in the loss function of the neural network.  
To allow the neural network to operate uniformly with respect to the small scales, it is desirable that the neural network satisfies an Asymptotic-Preservation (AP) property in the learning process. To this end, we consider a new class of AP Neural Networks (APNNs) for multiscale hyperbolic transport models of epidemic spread that, thanks to an appropriate AP formulation of the loss function, is capable to work uniformly at the different scales of the system. A series of numerical tests for different epidemic scenarios confirms the validity of the proposed approach, highlighting the importance of the AP property in the neural network when dealing with multiscale problems especially in presence of sparse and partially observed systems. 
\end{abstract}


\tableofcontents

\section{Introduction}
In recent decades, enormous progress has been made in the understanding of complex systems described by multiscale PDEs with applications ranging from classical physics and engineering to biology and social sciences^^>\cite{ABBDPTZ21,PareschiToscani,Bellomo,Bellomo2,BDP,Jin2022,Carrillo}.
 
Despite continuing progress, modeling and predicting the evolution of nonlinear multiscale systems using classical analytical or computational tools inevitably faces severe challenges. Firstly, numerically solving a multiscale problem requires complex and sophisticated computational codes and can introduce prohibitive costs (due to the well-known \textit{curse of dimensionality}). Moreover, we are always facing the difficulties related to the scarcity of data and multiple sources of uncertainty, especially when concerning social sciences. Above all, solving real physical problems with missing or incomplete initial or boundary conditions through traditional approaches is currently impractical. This is where and why data-driven models began to play a crucial role^^>\cite{E2021}.

Machine Learning (ML) is an incredibly powerful tool, which has proven to have an enormous impact in many fields of our society. This has led to great interest in using ML techniques also to study challenging scientific problems in science, engineering, medicine, concerning complex multiscale dynamics.
However, it is clear that the problems we are dealing with are very different from the classical problems in which ML has proved to be so successful. So, we cannot simply take the available ML methods as a ``black box" and use them uncritically^^>\cite{Willcox}.
Purely data-driven models may fit observations very well, but predictions may result physically inconsistent and, consequently, lead to erroneous generalizations. Therefore, there is an urgent need to integrate fundamental physical laws and related mathematical models into  the learning process of the neural networks^^>\cite{Karniadakis2021,E2020,Lou2021}. 
The main motivation for developing this new class of physics-informed machine learning algorithms is that such prior physical knowledge or constraints can ensure that ML methods remain robust even in the presence of imperfect data (such as missing, incomplete or noisy data) and provide accurate predictions that adhere to the physics of the phenomenon under study.

A recent example of  this new learning paradigm is represented by Physics-Informed Neural Networks (PINNs)^^>\cite{Coutinho2022,Raissi2019,Karniadakis2021,Yu2021}. PINNs are a new class of deep neural networks (DNNs) that are trained to solve supervised learning tasks while respecting any given physical laws described through general nonlinear ordinary differential equations (ODEs) or partial differential equations (PDEs). The physical knowledge of the underlying phenomenon  is incorporated into the PINN mainly in two ways: either it is introduced directly through the data embodying the underlying physics of the phenomenon of interest (observational bias) or it is introduced by an appropriate choice of the loss function that the PINN must minimize, forcing the training phase of the neural network to converge to solutions that adhere to the underlying physics (learning bias).

Nevertheless, the adoption of a standard formulation of PINNs in the context of multiscale problems may still lead to incorrect inferences and predictions^^>\cite{Jin2021}. This is mainly due to the presence of small scales leading to reduced or simplified models in the system that need to be enforced consistently during the learning process. In these cases, a standard PINN formulation allows an accurate description of the process only at the leading order, thus loosing accuracy in the asymptotic limit regimes. One remedy for this, as recently proposed in^^>\cite{Jin2021}, is to modify the loss function to include asymptotic-preserving (AP) properties during the training process. The realization of such an AP-loss function will therefore depend on the particular problem under study and will be based on an appropriate asymptotic analysis of the model.   

One particularly interesting area where the use of machine learning techniques can play a key role concerns epidemiological dynamics. In this context, a number of mathematical models have recently been proposed that require the estimation of several parameters from data to provide predictive scenarios and to test their reliability^^>\cite{ABBDPTZ21,Bert2,Bert3,BDP,Buonomo,Gatto,Giordano,Parolini,Veneziani2021}.
In this paper we will focus on a new class of epidemic models described by multiscale PDEs capable of describing both hyperbolic-type phenomena characteristic of epidemic propagation over long distances and main lines of communication between cities and parabolic-type phenomena in which classical diffusion prevails at the urban level^^>\cite{ABBDPTZ21,Bert2,Bert3,BDP}. The multiscale nature of the problem poses a challenge to the construction of PINN, and preservation of the AP property is therefore essential in order to obtain reliable results. Following the approach recently introduced in^^>\cite{Jin2021}, we will show how to construct AP neural networks (APNNs) that are capable to solve both inverse and forward problems of interest in epidemic dynamics.

The rest of the paper is organized as follows. The next section is devoted to the description of the model under study and a formal analysis of the different multiscale behaviors. In Section 3 we introduce the notion of APNN and describe how to construct such a neural network in the case of a simplified multiscale hyperbolic model and then how to extend it to the epidemic case under study. A series of numerical results for both inverse and forward problems using synthetic data produced by the numerical solution of the mathematical model illustrate the validity of the present approach. In particular, the case of partially observed systems, as commonly found in epidemics, will be considered and permits to emphasize the relevance of the AP-property. Some final considerations and future developments are reported in a concluding section.

\section{Hyperbolic models of epidemic spread}
\label{section_SIR}
For simplicity, we illustrate the space dependent epidemiological modeling in the case of a classic SIR compartmental dynamics, in which we consider a population subdivided in susceptible $S$ (individuals who may be infected by the disease), infectious $I$ (individuals who may transmit the disease) and removed $R$ (individuals healed and immune or died due to the disease). We assume to have a population with subjects having no prior immunity and neglect the vital dynamics represented by births and deaths due to the time scale considered. Nevertheless, it is straightforward to extend our arguments to more enriched compartmentalizations, designed to take into account specific features of the infectious disease of interest, as those proposed recently in^^>\cite{Bert2,Bert3,BDP,Buonomo,Gatto,Giordano,Parolini,Veneziani2021} to study the spread of COVID-19. 

\subsection{The hyperbolic SIR model}
By analogy with discrete-velocity kinetic theory^^>\cite{PareschiToscani,Bellomo}, we consider individuals moving in a one-dimensional domain $\mathcal{D} \subseteq \mathbb{R}$ in two opposite directions, with velocities $\pm \lambda_{S,I,R} = \pm \lambda_{S,I,R}(x)$, distinguished for each epidemic compartment. Notice that the characteristic velocities reflect the heterogeneity of geographical areas, and, therefore, are chosen dependent on the spatial location $x \in \mathcal{D}$. Hence, we can describe the space-time dynamics of the population for $t>0$ through the following two-velocity SIR epidemic transport model^^>\cite{ABBDPTZ21,Bert,Bert4}:
\begin{equation}
\begin{split}
	\frac{\partial S^{\pm}}{\partial t} + \lambda_S \frac{\partial S^{\pm}}{\partial x} &= - \beta S^{\pm}I \mp \frac{1}{2\tau_S}\left(S^+ - S^-\right),	\\ 
	\frac{\partial I^{\pm}}{\partial t} + \lambda_I \frac{\partial I^{\pm}}{\partial x} &=  \beta S^{\pm}I -\gamma I^{\pm} \mp \frac{1}{2\tau_I}\left(I^+ - I^-\right),\\ 
	\frac{\partial R^{\pm}}{\partial t} + \lambda_R \frac{\partial R^{\pm}}{\partial x} &= \gamma I^{\pm} \mp \frac{1}{2\tau_R}\left(R^+ - R^-\right),				
\end{split}
\label{eq.SIR_kinetic_diag}
\end{equation}
with the total densities of each compartment, $S(x,t)$, $I(x,t)$, and $R(x,t)$, given by
\begin{equation*}
 	S = S^+ + S^- , \quad 
 	I = I^+ + I^- , \quad 
 	R = R^+ + R^- .	
\end{equation*}
The transport dynamics of the population is governed by the scaling parameters $\lambda_{S,I,R}$ as well as the relaxation times $\tau_{S,I,R}=\tau_{S,I,R}(x)$. The quantity $\gamma=\gamma(x,t)$ is the recovery rate of infected, which corresponds to the inverse of the infectious period. This rate may vary in space and time depending on the treatment therapies used, even though generally can be assumed constant, especially for short-term analysis. The transmission of the infection is defined by an incidence function $\beta S I$ modeling the transmission of the disease^^>\cite{HWH00,CS78,Kermack}.
The transmission rate $\beta=\beta(x,t)$ characterizes the average number of contacts per person per time, multiplied by the probability of disease transmission in a contact between a susceptible and an infectious subject. Notice that  this rate may vary in space and time as a consequence of the intensification of governmental control actions (such as mandatory wearing of masks, closing of specific activities or full lockdowns) or their relaxation (lifting mask mandates, reopening schools, restaurants, leisure and cultural centers) in specific locations. 

It is worth to highlight that, when investigating real epidemic scenarios, the  above-mentioned  parameters are, in general, unknown. While the recovery rate might be fixed based on clinical data, the transmission rate must always be estimated through a delicate calibration process in order to match available data. It is also well-known that this process is highly heterogeneous which makes the inverse problem even more challenging \cite{Dimarco2021}.

The standard threshold of epidemic models is the so-called basic reproduction number $R_0$, which defines the average number of secondary infections produced by one infected individual in a totally susceptible population^^>\cite{HWH00}. The effective reproduction number $R_t$, instead, defines the variation in time of this rate, giving information on the progress of the infectious spread. Indeed, this number determines when an infection can invade and persist in a new host population ($R_t > 1$), or tend to fade away ($R_t < 1$). The endemic state corresponds to the case $R_t = 1$.

Assuming no inflow/outflow boundary conditions in $\mathcal{D}$, integrating in space and summing up the second equation in \eqref{eq.SIR_kinetic_diag} we are able to define the effective reproduction number of the SIR transport model
\begin{equation}
R_t(t)=\frac{\int_{\mathcal{D}} \beta(x,t) S(x,t) I(x,t)\,dx}{\int_{\mathcal{D}} \gamma(x,t) I(x,t)\,dx} \geq 1.
\label{eq:R0}
\end{equation}
Notice that this definition naturally extends locally by integrating over any subset of the computational domain $\mathcal{D}$ if one ignores the boundary flows. Under the same no inflow/outflow boundary conditions, if we integrate in $\mathcal{D}$ equations \eqref{eq.SIR_kinetic_diag}, we can finally observe that the model fulfill the conservation of the total population, being
\begin{equation}
\frac{\partial}{\partial t} \int_{\mathcal{D}} \left(S(x,t)+I(x,t)+R(x,t)\right)\,dx =0\,, 
\label{eq:conservation}
\end{equation}
with $S(t)+I(t)+R(t)=P$, and $P$ total population reference size, constant over time.

\subsection{Multiscale behavior and diffusion limit}
\label{sect:diff-LF}
Introducing the fluxes, defined by
\begin{equation}
 	J_S = \lambda_S \left(S^+ - S^-\right) , \quad 
 	J_I = \lambda_I \left(I^+ - I^-\right) , \quad 
 	J_R = \lambda_R \left(R^+ - R^-\right) ,	  	\label{eq.fluxes_kinetic_diag}
\end{equation}
we obtain a hyperbolic model equivalent to \eqref{eq.SIR_kinetic_diag}, but presenting a macroscopic description of the propagation of the epidemic at finite speeds 
\begin{equation}
\begin{split}
	\frac{\partial S}{\partial t} + \frac{\partial J_S}{\partial x} &= -\beta SI  ,			\\ 
	\frac{\partial I}{\partial t} + \frac{\partial J_I}{\partial x} &= \beta SI  -\gamma I,	\\ 
	\frac{\partial R}{\partial t} + \frac{\partial J_R}{\partial x} &= \gamma I , \\
	\frac{\partial J_S}{\partial t} + \lambda_S^2 \frac{\partial S}{\partial x} &= -\beta J_SI  -\frac{J_S}{\tau_S},\\ 
	\frac{\partial J_I}{\partial t} + \lambda_I^2 \frac{\partial I}{\partial x} &= \frac{\lambda_I}{\lambda_S}\beta J_SI  -\gamma J_I -\frac{J_I}{\tau_I},\\ 
	\frac{\partial J_R}{\partial t} + \lambda_R^2 \frac{\partial R}{\partial x} &= \frac{\lambda_R}{\lambda_I}\gamma J_I -\frac{J_R}{\tau_R} .
\end{split}
\label{eq.SIR_kinetic}
\end{equation}
Let us now consider the behavior of this model in diffusive regimes^^>\cite{LT}. To this aim, we introduce the space dependent diffusion coefficients 
\begin{equation}
D_S=\lambda_S^2 \tau_S,\quad D_I=\lambda_I^2 \tau_I,\quad D_R=\lambda_R^2 \tau_R.
\label{eq:diff2}
\end{equation}
which characterize the diffusive transport mechanism of susceptible, infectious and removed, respectively. Keeping the above quantities fixed while letting the relaxation times $\tau_{S,I,R}\to 0$ (and so the characteristic velocities $\lambda_{S,I,R} \to \infty$), from the last three equations in \eqref{eq.SIR_kinetic} we obtain, for each epidemic compartment, a proportionality relation between the flux and the spatial derivative of the corresponding density (Fick's law)
\begin{equation}
J_S = -D_S \frac{\partial S}{\partial x},\quad  J_I = -D_I\frac{\partial I}{\partial x},\quad J_R = -D_R\frac{\partial R}{\partial x}.
\label{eq:fick}
\end{equation}
Substituting \eqref{eq:fick} into the first three equations in \eqref{eq.SIR_kinetic}, we recover the following parabolic reaction-diffusion model, widely used in literature to study the spread of infectious diseases^^>\cite{MWW,Sun,Webb,HS,Berestycki}
\begin{equation}
\begin{split}
\frac{\partial S}{\partial t} &=  -\beta SI +\frac{\partial}{\partial x} \left({D_S}\frac{\partial S}{\partial x}\right),\\
\frac{\partial I}{\partial t} &=  \beta SI-\gamma I+\frac{\partial}{\partial x} \left({D_I}\frac{\partial I}{\partial x}\right),\\
\frac{\partial R}{\partial t} &=  \gamma I+\frac{\partial}{\partial x}\left({D_R}\frac{\partial R}{\partial x}\right).
\end{split}
\label{eq:diff}
\end{equation}

The model's capability to account for different regimes, ranging from hyperbolic to parabolic, according to the space dependent values $\tau_{S,I,R}$ and $\lambda_{S,I,R}$, makes it suitable for describing the dynamics of human beings. Our daily routine is, indeed, a complex mixing of individuals moving at the scale of a city center and individuals traveling among different municipalities. 
In this situation, it results more appropriate to describe the human dynamics in city centers with a high density of individuals through a diffusion operator, while characterizing the mobility of subjects in extra-urban areas through a hyperbolic, advective, mechanism, avoiding in this case a propagation of the information at infinite speeds^^>\cite{Bert2,Bert3,BDP}. 

\section{Asymptotic-Preserving Neural Networks (APNNs)}
In this section, we provide a brief overview of the general framework of PINNs \cite{Karniadakis2021,Raissi2019} and then we shall discuss the relevant concepts of  {\it Asymptotic-Preserving Neural Networks} (APNNs) for the problems of interest.

\subsection{Basics of PINNs}
The design of a standard deep neural network (DNN) by supervised learning can be summarized in three main steps^^>\cite{E2020}: 
\begin{enumerate}
\item The choice of the neural network structure.  
\item The loss function that minimizes the classical empirical risk, typically characterized by the difference between model and data.
\item A method to minimize loss over the parameter space. The most popular choices are stochastic gradient descent (SGD) and advanced optimizers such as Adam^^>\cite{Adam}.
\end{enumerate}
In practice, the performance of the neural network is estimated on a finite data set (which is unrelated to any data used to train the model) and called \emph{test error}, whereas the error in the loss function (which is used for training purposes) is called the \emph{training error}.

Compared to the above classical deep learning methodology, the major difference of PINN is the integration of physical laws, usually in the form of PDEs
\begin{equation}
\label{eq:general-pde-system}
\begin{aligned}
    \mathcal{F}(u, x, t; \xi) &= 0, \quad && (x, t) \in \Omega, \\
    \mathcal{B}(u, x, t; \xi) &= 0, \quad && (x, t) \in \partial \Omega.
\end{aligned}
\end{equation}
Here $\Omega \subset \mathbb{R}^{d}\times\mathbb{R}$ is the spatio-temporal domain of the system, $\partial \Omega$ represents its boundary, $\mathcal{F}$ is the differential operator, $u$ represents the solution to the system, $\xi$ is the parameter related to the physics. Since the initial condition is mathematically equivalent to the boundary condition in the spatio-temporal domain, we use $\mathcal{B}$ as a general operator for arbitrary initial and boundary conditions of the system. 

PINN models usually include a neural network representation of the solution $u \approx u_{NN}(z; \theta)$, parameterized by network parameters $\theta$ and having $z\in\mathbb{R}^d$ input data. In PINN literature, the most widely used neural network architecture is the feed-forward neural network (FNN). 
A $L+1$ layered FNN consists of an input layer, an output layer, and $L-1$ hidden layers, which can be defined as follows
\[
\begin{split}
z^1 &= W^1 z + b^1,\\
z^l &= \sigma \circ (W^l z^{l-1} + b^l),\quad l=2,\ldots,L-1\\
u_{NN}(z;\theta) & = z^L= W^L z^{L-1} + b^L,
\end{split}
\]
where $W^{l}\in\mathbb{R}^{m_{l}\times m_{l-1}}$ are the weights, $b^l\in\mathbb{R}^{m_{l}}$ the bias, $m_l$ is the width of the $l$-th hidden layer with $m_1 = d_{in} = d$ the input dimension and $m_L = d_{out}$ the output dimension,  $\sigma$ is a scalar activation function (such as ReLU \cite{Goodfellow-et-al-2016}), and ``$\circ$'' denotes entry-wise operation. Thus, we denote the set of network parameters $\theta=(W^1,b^1,\ldots,W^L,b^L)$.

To find the optimal values for $\theta$, the neural network is trained by minimizing the following type of loss (also called cost or risk) function
\begin{equation}
\label{eq:general-pinn-loss}
\mathcal{L}(\theta) = w_r^T \mathcal{L}_{r}(\theta) + w_b^T \mathcal{L}_{b}(\theta) +  w_d^T \mathcal{L}_{d}(\theta).
\end{equation}
Here $\mathcal{L}_{r}$ and $\mathcal{L}_{b}$ quantify the discrepancy of the neural network surrogate $u_{NN}$ with the underlying PDE and its initial or boundary conditions in^^>\eqref{eq:general-pde-system}, respectively. The data mismatch loss $\mathcal{L}_{d}$ is applied when additional measurement data are available, e.g., when solving inverse problems, and $w_{r}$, $w_{b}$, $w_{d}$ are the corresponding weight vectors. The most popular methods chosen to solve this optimization problem remain stochastic gradient descent (SGD) and Adam^^>\cite{Adam}. 
After finding the optimal set of parameter values $\theta^*$ by minimizing the PINN loss \eqref{eq:general-pinn-loss}, i.e.,
\begin{equation}
\theta^* = \argmin \mathcal{L}(\theta),
\end{equation}
the neural network surrogate $u_{NN}( x,t; \theta^*)$ can be evaluated at any given spatio-temporal point to get the solution. 

In the context of inverse problems, the structure of the network is almost the same with respect to the forward problem setting, except that unknown physical parameters $\xi$ are treated as learnable parameters. As a result, the training process involves optimizing $\theta$ and $\xi$ jointly
\begin{equation}
    (\theta^*, \xi^*) = \argmin \mathcal{L}(\theta,\xi).
\end{equation}
In summary, PINN can be regarded as an unsupervised learning approach when used to train forward problems, with only equations residual and boundary conditions in the loss function, and as a semi-supervised learning  approach for inverse problems, when some measurements are available. In the last part of this section, we shall further discuss in detail each component of this learning framework  through several examples. 

\subsection{Extension to APNNs}

\begin{figure}[t!]
\centering
    \includegraphics[width=0.6\textwidth]{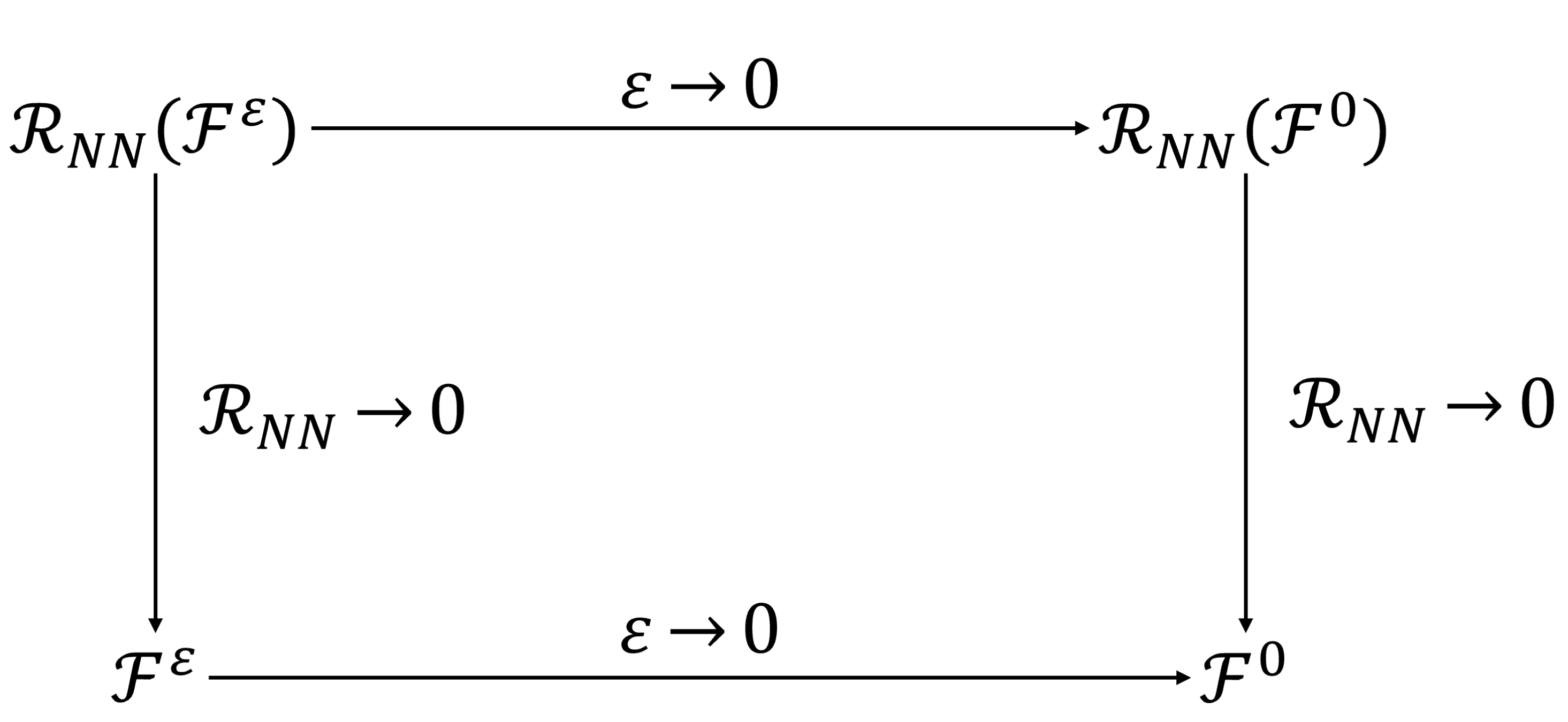}
    \caption{AP diagram for neural networks. $\mathcal{F}^\epsilon$ is the multiscale problem that depends on the scaling parameter $\epsilon$, while $\mathcal{F}^0$ is the corresponding formulation in the reduced order limit, which does not depend anymore on $\epsilon$. The solution of the system $\mathcal{F}^\epsilon$ is approximated by the neural network through the imposition of the residual term $\mathcal{R}_{NN}\left( \mathcal{F}^\epsilon\right) = \mathcal{R}_{NN}^\epsilon$. The asymptotic limit of $\mathcal{R}_{NN}\left( \mathcal{F}^\epsilon\right)$ as $\epsilon \to 0$ is denoted with $\mathcal{R}_{NN}\left( \mathcal{F}^0\right)  = \mathcal{R}_{NN}^0$. The neural network is called AP if $\mathcal{R}_{NN}\left( \mathcal{F}^0\right)$ is consistent with the residual of the reduced system $\mathcal{F}^0$.}
    \label{fig:AP_scheme}
\end{figure}

Since we aim at analyzing multiscale hyperbolic dynamics regardless of the propagation scaling, in order to obtain physically-based predictions, it is important that the PINN can preserve the correct equilibrium solution \eqref{eq:diff} in the diffusive regime, which means that the PINN should fulfill the AP property^^>\cite{Jin2022,Jin2021,Lu2021,Han2019}. We remark that in the context of the epidemic modeling of this work, the AP property is of particular importance, allowing the same neural network to efficiently and robustly simulate population dynamics characterized by both diffusive and hyperbolic transport behaviors (the former in urban centers and the latter for mobility along connecting routes). 

The neural networks satisfying this property are called {\it Asymptotic-Preserving Neural Networks} (APNNs), and have been recently introduced in \cite{Jin2022,Jin2021} to efficiently solve multiscale kinetic problems with scaling parameters that can have several orders of magnitude of difference. The definition of an APNN reported in \cite{Jin2021} for the case of multiscale kinetic models with continuous velocity fields is generalized in the following (see Figure \ref{fig:AP_scheme}).

\begin{definition}[Asymptotic-Preserving Neural Network]
Assume the solution is parameterized by a PINN trained by using an optimization method to minimize a loss function which includes a residual term enforcing the physics of the phenomenon. Then we say it is an \textit{Asymptotic-Preserving Neural Network} (APNN) if, as the physical scaling parameter of the multiscale model tends to zero, the loss function of the full model-constraint converges to the loss function of the corresponding reduced order model. 
\end{definition}
In other words, the loss function, viewed as a numerical approximation of the original equation, benefits from the AP property.

\subsection{A simple example: APNN for the Goldstein--Taylor model}
\label{sec:test1-goldstein-taylor}
To illustrate the relevance of the AP property in the construction of the neural network, let us carry on a detailed example by considering a simplified case in which there are no epidemic source terms that allow individuals to move to a different compartment and the entire population behaves as a single compartment. Such a case corresponds to the so-called Goldstein--Taylor model in discrete velocity kinetic theory^^>\cite{PareschiToscani,JPT}. This model, indeed, describes the space-time evolution of the two particles densities $f^+(x,t)$ and $f^-(x,t)$, at time $t>0$, traveling in a one-dimensional domain, $x \in \mathcal{D} \subseteq \mathbb{R}$, with velocity $\pm c$, respectively. At the same time, particles can change and assume the opposite velocity, randomly. The dynamics of this system of particles is governed by the following system of PDEs
\begin{equation}
\begin{split}
\frac{\partial f^\pm}{\partial t} \pm \frac{c}{\epsilon} \frac{\partial f^\pm}{\partial x} &= \frac{\sigma}{2\epsilon^2} \left(f^\mp - f^\pm\right),	
\end{split}
\label{Goldstein-Taylor}
\end{equation}
with $\epsilon$ scaling parameter of the kinetic dynamics and $\sigma$ scattering coefficient. The total particles density is given by $\rho(x,t)=f^+(x,t)+f^-(x,t)$. 

We consider ${f}_{NN}^\pm(x,t;\theta)$ to be a DNN with inputs $x$ and $t$ and trainable parameters $\theta$, to  approximate the solution of our system: ${f}^\pm(x,t) \approx {f}_{NN}^\pm(x,t;\theta)$. Then, we define the PDEs residual 
\begin{equation}
    \label{eq:pinn-nonap-residue1}
    \mathcal{R}^\epsilon_{NN}(f^\pm) = \epsilon^2\frac{\partial f^\pm_{NN}}{\partial t} \pm \epsilon c \frac{\partial f^\pm_{NN}}{\partial x} - \frac{\sigma}{2} \left(f^\mp_{NN}-f^{\pm}_{NN}\right),
\end{equation}
and incorporate it into the loss function term $\omega_r^T {\mathcal L}_r(\theta)$ of the neural network by taking the weighted mean square error of the residual to obtain a standard PINN. 

To understand the asymptotic behavior of the model we resort on a suitable macroscopic formulation of the system which is achieved through the introduction of the scaled flux $j = c\left(f^+ - f^- \right)/\epsilon$. This permits to write the system \eqref{Goldstein-Taylor} in equivalent form as
\begin{equation}
\begin{split}
\frac{\partial \rho}{\partial t} +  \frac{\partial j}{\partial x} &= 0\,,	\\ 
\frac{\partial j}{\partial t} + \frac{c^2}{\epsilon^2} \frac{\partial \rho}{\partial x} &= -\frac{\sigma}{\epsilon^2} j \,.
\end{split}
\label{Goldstein-Taylor_macro}
\end{equation}
In the diffusion limit, i.e. let $\epsilon \to 0$, we obtain 
\begin{equation}
\begin{split}
j &= -\frac{c^2}{\sigma}\frac{\partial \rho}{\partial x},
\end{split}
\label{eq.j}
\end{equation}
which, inserted into the first equation, leads to the reduced diffusive model (which recalls the standard heat equation)
\begin{equation}
\begin{split}
\frac{\partial \rho}{\partial t} &= \frac{c^2}{\sigma} \frac{\partial^2 \rho }{\partial x^2}.
\end{split}
\label{eq.heat}
\end{equation}
It is clear that the standard PINN residual \eqref{eq:pinn-nonap-residue1} is not consistent with the above analysis since $\mathcal{R}^\epsilon_{NN}(f^\pm)$ in the limit $\epsilon\to 0$ reduces to 
\[
\mathcal{R}^0_{NN}(f^\pm ) = - \frac{\sigma}{2} \left(f^\mp_{NN}-f^{\pm}_{NN}\right),
\]
which corresponds to force $f^+(x,t)=f^-(x,t)$ and does not suffice to achieve the correct diffusive behavior \eqref{eq.heat}.

In contrast, using the macroscopic formulation \eqref{Goldstein-Taylor_macro}, we can construct an APNN incorporating in the loss function the mean square error of the PDEs residuals
\begin{equation}
    \label{eq:pinn-ap-residue_GT}
    \begin{split}
    \mathcal{R}^\varepsilon_{NN}(\rho) =  \frac{\partial \rho_{NN}}{\partial t} +  \frac{\partial j_{NN}}{\partial x},\qquad
    \mathcal{R}^\varepsilon_{NN}(j) =  \varepsilon^2 \frac{\partial j_{NN}}{\partial t} +  c^2 \frac{\partial \rho_{NN}}{\partial x}+\sigma j_{NN}.
    \end{split}
\end{equation}
Now, in the limit $\epsilon \to 0$, we obtain
\begin{equation}
    \label{eq:pinn-ap-residue_GTl}
    \begin{split}
    \mathcal{R}^0_{NN}(\rho) =  \frac{\partial \rho_{NN}}{\partial t} +  \frac{\partial j_{NN}}{\partial x},\qquad
    \mathcal{R}^0_{NN}(j) =  c^2 \frac{\partial \rho_{NN}}{\partial x}+\sigma j_{NN},
    \end{split}
\end{equation}
which is consistent with the residual of the limiting diffusive model \eqref{eq.heat}.
We refer to Appendix \ref{appendix} for a detailed description of the loss function for the Goldstein-Taylor model, including data and boundary conditions loss terms.

\subsection{APNN for the hyperbolic SIR model}
To achieve the AP property in the neural network for the hyperbolic SIR model, we follow the same approach of the previous section. Thus, we  consider the system written in macroscopic form defined by equations \eqref{eq.SIR_kinetic}. Multiplying both members of each equation for the corresponding scaling parameter $\tau_i, i\in \{S,I,R\}$, we can rewrite the system in the following compact form
\begin{equation}
    \label{eq:ap-vector}
    \begin{aligned}
    {\mathbf{\tau}}(x)\frac{\partial {U}(x,t)}{\partial t} + {D}(x) \frac{\partial {F}({U(x,t)})}{\partial x} &= {G}({U(x,t)}), \quad && (x, t)\in \Omega, 
    \end{aligned}
\end{equation}
where 
\begin{equation*}
    \label{eq:ap-vector-components}
    {U} = 
        \begin{bmatrix}
            S \\ I \\ R \\ J_S \\ J_I \\ J_R
        \end{bmatrix}
    , \
    {\tau} = 
        \begin{bmatrix}
            1 \\ 1 \\ 1 \\ \tau_S \\ \tau_I \\ \tau_R
        \end{bmatrix}
    , \
    {D}= 
        \begin{bmatrix}
            1 \\ 1 \\ 1 \\ D_S \\ D_I \\ D_R
        \end{bmatrix}, \
    {F}({U}) = \begin{bmatrix}
            J_S \\ J_I \\ J_R \\ S \\ I \\ R
    \end{bmatrix}, \
    {G}({U}) = \begin{bmatrix}
            -\beta SI \\
            \beta SI -\gamma I \\
            \gamma I \\
            -\tau_S\beta J_S I - J_S \\
            \tau_I\frac{\lambda_I}{\lambda_S}\beta J_S I - \tau_I\gamma J_I - J_I \\
            \tau_R\frac{\lambda_R}{\lambda_I}\gamma J_I - J_R
        \end{bmatrix}.
\end{equation*}

We consider ${U}_{NN}(x,t;\theta)$ to be a deep neural network (NN) with inputs $x$ and $t$ and trainable parameters $\theta$, to  approximate the solution of our system: ${U}(x,t) \approx {U}_{NN}(x,t;\theta)$. Then, we define the residual term
\begin{equation}
    \label{eq:pinn-ap-residue}
    \mathcal{R}^{\tau}_{NN}(U) = {\tau}\frac{\partial {U}_{NN}}{\partial t} + {D} \frac{\partial {F}({U}_{NN})}{\partial x} - {G}({U}_{NN}),
\end{equation}
and embed it into the loss function of the neural network to obtain an APNN. We omit for brevity the detailed analysis of the AP property. In the limit as $\tau_i \to 0$, $\lambda_i\to \infty$, $i\in\{S,I,R\}$, under conditions \eqref{eq:diff2}, such analysis follows the same steps of the previous section, and $\mathcal{R}^{0}_{NN}$ results in agreement with the diffusion limit computed in Section \ref{sect:diff-LF}.
 
\begin{figure}[t!]
\centering
    \includegraphics[width=\textwidth]{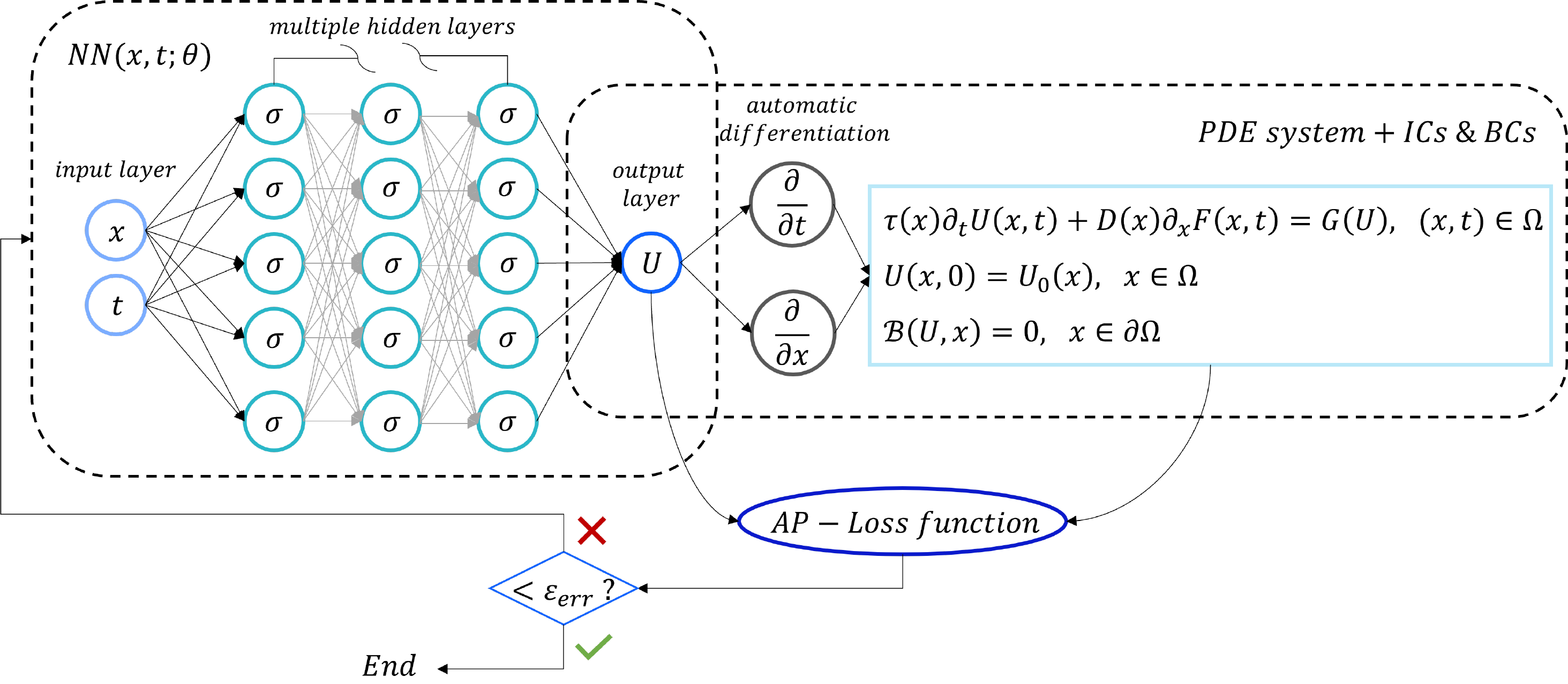}
    \caption{APNN schematic work-flow. The NN architecture is integrated with the physical knowledge of the dynamics of interest through the inclusion of the PDE system and the enforcement of initial and boundary conditions (and eventually conservation properties), when known, becoming a PINN. The AP property, which is a fundamental feature when dealing with multiscale hyperbolic systems, is guaranteed through the correct design of an \emph{AP-loss} function.}
    \label{fig:APNN_scheme}
\end{figure}
We restrict the neural network approximation ${ U}_{NN}$ to satisfy the physics imposed by the residual \eqref{eq:pinn-ap-residue} on a finite set of $N_r$ user-specified scattered points inside the domain, $\{(x_r^n, t_r^n)\}_{n=1}^{N_r} \subset \Omega$ (referred as \textit{residual points}) and we also enforce the initial and space-boundary conditions of the system on $N_b$ scattered points of the space-time boundary $\mathcal{B}({ U}(x, t))$, 
i.e. $\{(x_b^k, t_b^k)\}_{k=1}^{N_b} \subset \partial \Omega$ \cite{Kharazmi2021}.
In the context of inverse problems, we also consider to have access to measured data, with a dataset $\{({ U}_d^i, x_d^i, t_d^i)\}_{i=1}^{N_d}$, with ${ U}_d^i= { U}(x_d^i, t_d^i)$, available in a finite set of fixed \textit{training points}.
Thus, in the training process of the PINN, we minimize the following AP-loss function, composed of four mean squared error terms
\begin{equation}
    \label{eq:pinn-loss}
    \begin{aligned}
    \mathcal{L}(\theta) &= \omega_d^T \,{\mathcal L}_d(\theta) + \omega_b^T\, {\mathcal L}_b(\theta) + \omega_r^T \,{\mathcal L}^\tau_r(\theta) + \omega_c \, {\mathcal L}_c(\theta), 
    \end{aligned}
\end{equation}
where $\omega_d$, $\omega_r$, $\omega_b$, $\omega_c$ characterize the weights associated to each contribution. 
Notice that ${\mathcal L}_d$ quantifies the mismatch of the approximated solution with respect to known data samples, while ${\mathcal L}_b$, ${\mathcal L}^\tau_r$ and ${\mathcal L}_c$ represent the discrepancy in initial/boundary conditions of \eqref{eq:ap-vector}, in the residual \eqref{eq:pinn-ap-residue} and with respect to the conservation of the total density in the domain \eqref{eq:conservation}, respectively, all three contributing to enforce the physical structure of the problem. 
We present the detailed expression of each term in \eqref{eq:pinn-loss} in Appendix \ref{appendix1}. A schematic representation of the APNN architecture is given in Figure \ref{fig:APNN_scheme}.

\section{Numerical examples and applications}
In this section, various numerical tests are presented to assess the performance of the proposed APNNs. The first two examples concern the usage of an APNN for the solution of inverse and forward problem set up considering as prototype multiscale hyperbolic system either the standard Goldstein--Taylor model or a slightly modified version of it. Even if  this model is a simpler system of equations with respect to \eqref{eq.SIR_kinetic}, it well represents the dynamics of interest, as discussed in Section \ref{sec:test1-goldstein-taylor}. These tests are designed to further highlight how the choice of the APNN formulation proposed in this work is fundamental for the treatment of multiscale problems, especially in the context of availability of partial information. We shall demonstrate also numerically with this prototype model (and we refer to Section \ref{sec:test1-goldstein-taylor} for the analytical proof) that a standard PINN formulation leads to the loss of the AP property and, consequently, to non-physical reconstructions of the sought dynamics.

Following that, various tests concerning the solution of epidemic problems are discussed, examining the APNN performance in inferring the unknown epidemic parameters, solving the forward problem, and forecasting the spread of the infectious disease, also when spatially heterogeneous parameters are considered.

The numerical solution obtained with a second-order AP-IMEX Runge-Kutta Finite Volume method^^>\cite{Bert,Bert2} is considered as synthetic data for the ground truth and used in the APNN to build up the training dataset. With regards of epidemic test cases, we remark here, as also discussed in Appendix \ref{appendix1}, that since data of fluxes $J_S, J_I, J_R$ are not accessible in real-world applications, we only enforce the measurements of $S, I, R$ in ${\mathcal L}_d$. Nevertheless, unless otherwise specified, we impose initial conditions of the fluxes in ${\mathcal L}_b$. In all the examples, periodic boundary conditions are considered. To strictly impose them (accounted again in ${\mathcal L}_b$), we employ the periodic mapping technique taken from^^>\cite{zhang2020learning} in the input layer
\begin{equation}
    \label{eq:pinn-periodic-mapping}
    { U}_{NN}(x, t) = { U}_{NN}\left(\cos(\alpha x), \sin(\alpha x), t\right),
\end{equation}
where $\alpha$ is a hyperparameter controlling the frequency of the solution. 
For the tests concerning the Goldstein--Taylor model, the activation function \texttt{sin} is chosen, adopting the SIREN framework^^>\cite{sitzmann2020implicit}; for the epidemic tests, the function \texttt{tanh} is used. Finally, the Adam method^^>\cite{Adam} is used for the optimization process and derivatives in the NN are computed applying automatic differentiation^^>\cite{Baydin}.

For all the numerical examples, we adopt a single feed-forward neural network with depth $8$ and width $32$. The model structure is deliberately fixed among numerical experiments in both parabolic and hyperbolic regime, to highlight the main advantage of AP schemes that macroscopic behavior can be captured without resolving small physical parameters numerically (i.e. the architectural parameters of the neural network are independent of the physical scaling parameters). The chosen model and training hyperparameters are given in Tables \ref{Table:test12_hyperparameters} and \ref{Table:test34_hyperparameters} of the Appendix \ref{appendix_C} for each test case.


\subsection{Test 1: Goldstein-Taylor model in diffusive regimes}

In the following, we seek to emphasize numerically the importance of choosing the correct formulation to preserve the AP property and correctly approximate population dynamics even in diffusive regimes, particularly when dealing with partial information available. To this aim, we set up for problem \eqref{Goldstein-Taylor} a test with initial conditions
$$\rho(x,0) = 6 + 3\cos(3\pi x), \qquad j(x,0) = \frac{9\pi c^2}{\sigma}\sin(3\pi x),$$
with $c=1$ and $\sigma=4$. We consider periodic boundary conditions, choosing $\alpha = 3$ in the periodic mapping \eqref{eq:pinn-periodic-mapping}, and only the diffusive, parabolic regime of the model, choosing $\epsilon=10^{-4}$, with final time of the simulation $t_{end}=0.1$.

\begin{figure}[!t]
    \centering
    \includegraphics[width=.42\textwidth]{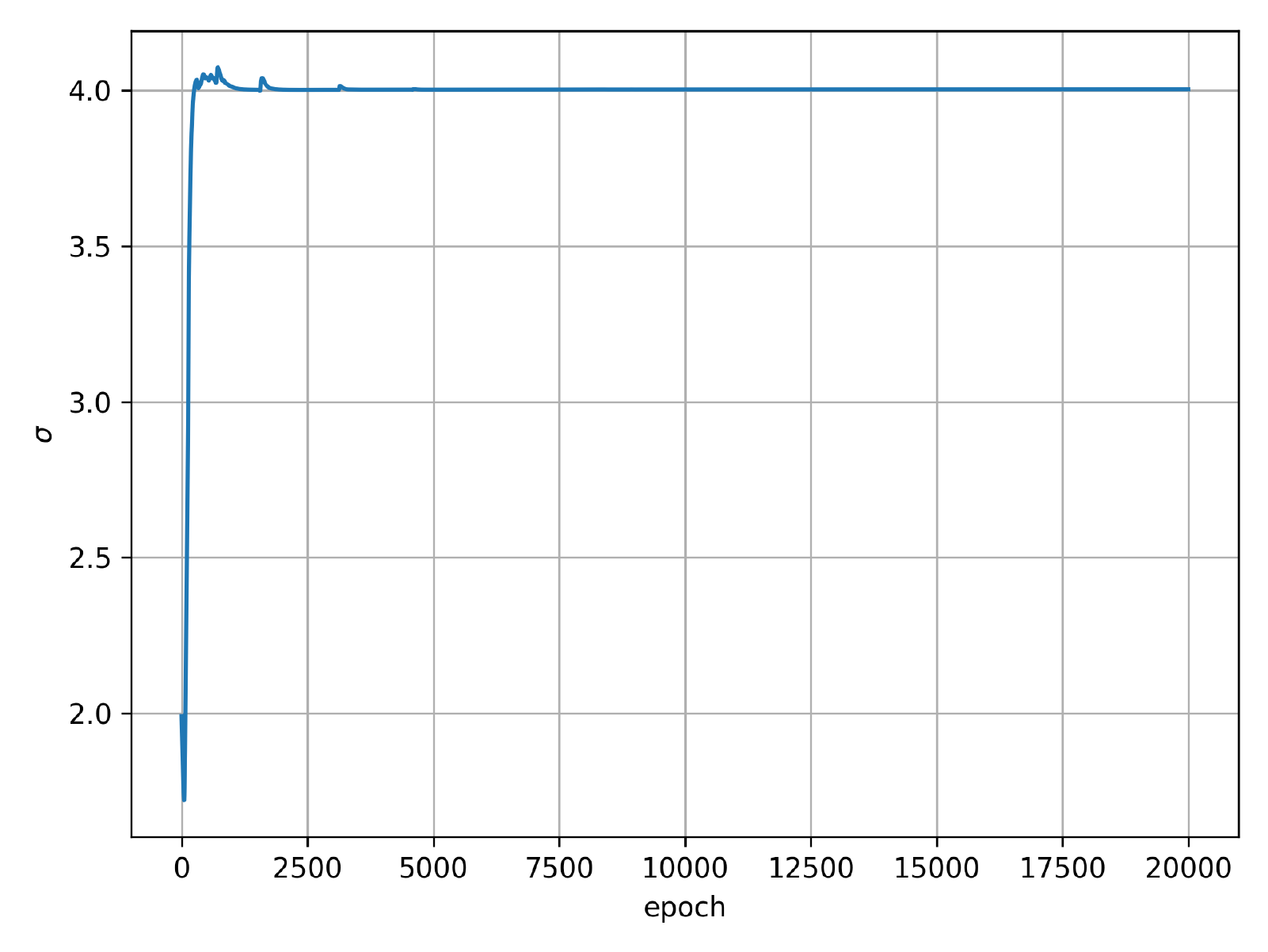}
    \includegraphics[width=.42\textwidth]{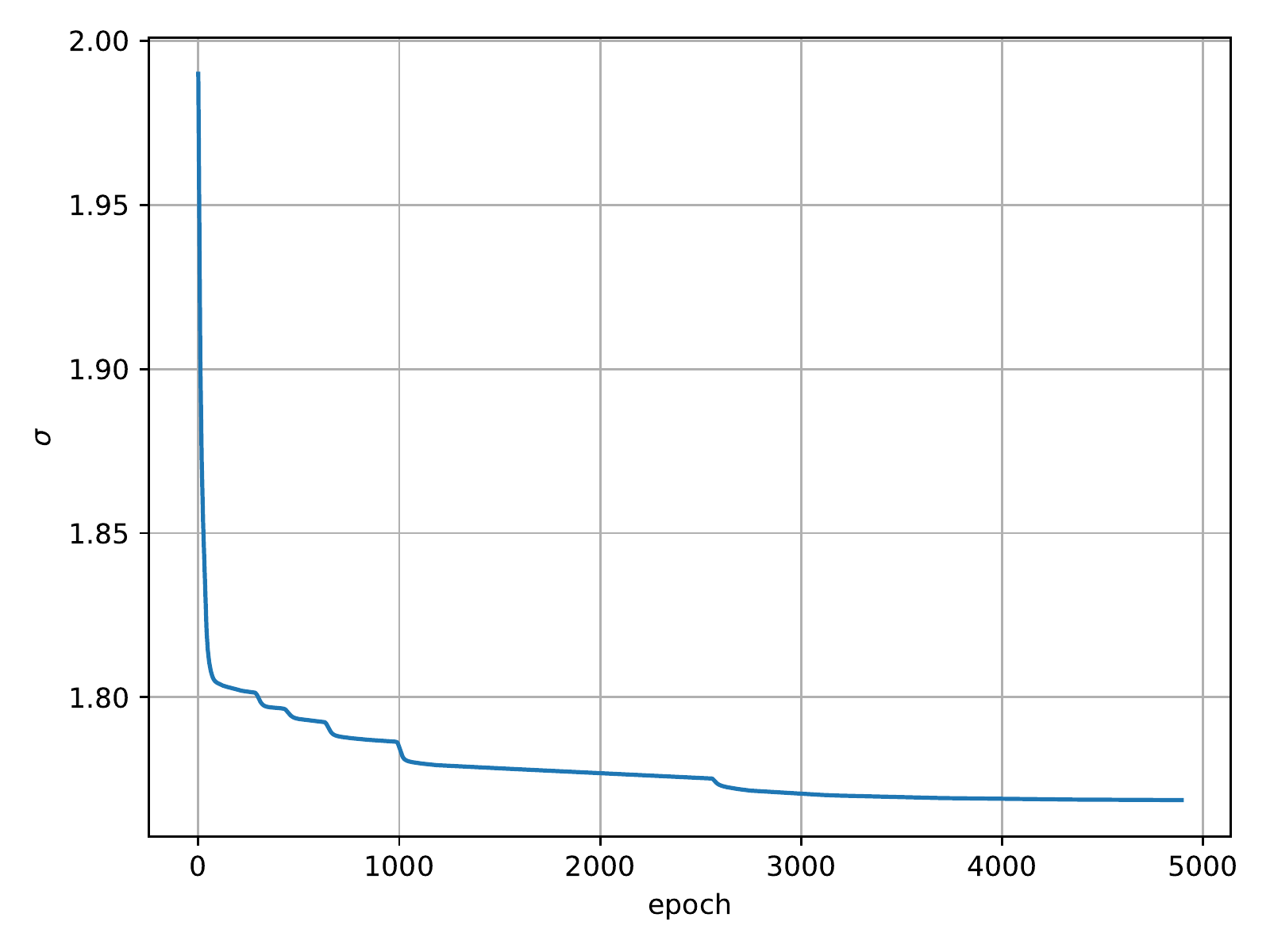}
    \caption{Test 1: Inverse problem for the Goldstein-Taylor model in the diffusive regime ($\epsilon=10^{-4}$). Convergence of the target parameter $\sigma=4$ with respect to epochs using the APNN (left) and the standard PINN (right).}
    \label{fig:test1_parabolic_sigma}
\end{figure}

\begin{figure}[!t]
    \centering
     \includegraphics[width=0.4\textwidth]{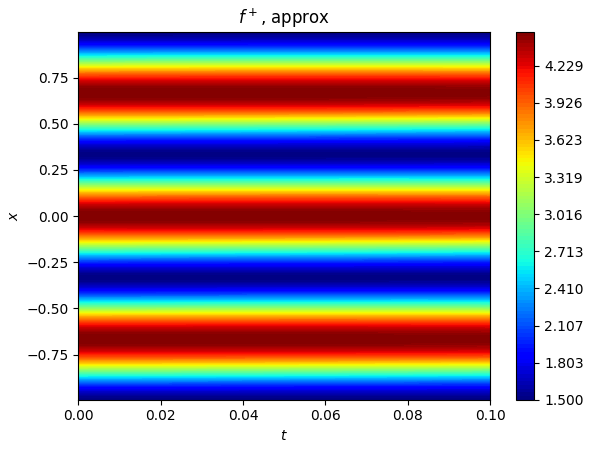}
     \includegraphics[width=0.4\textwidth]{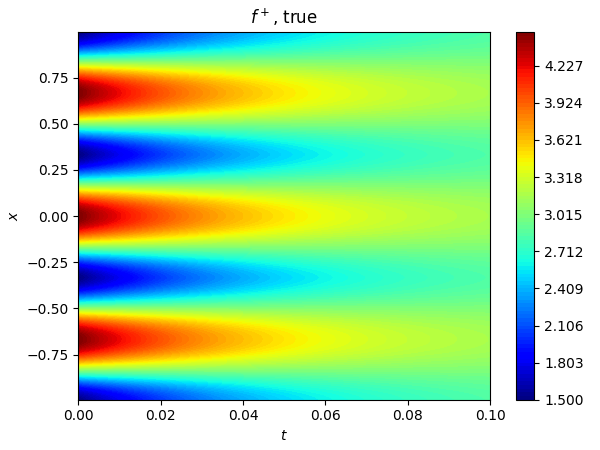}
    \includegraphics[width=0.4\textwidth]{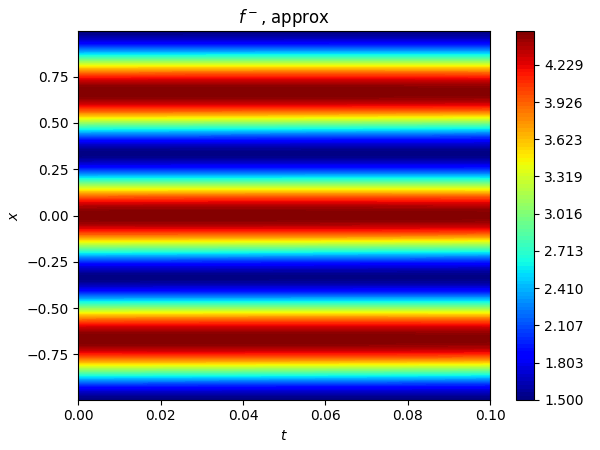}
     \includegraphics[width=0.4\textwidth]{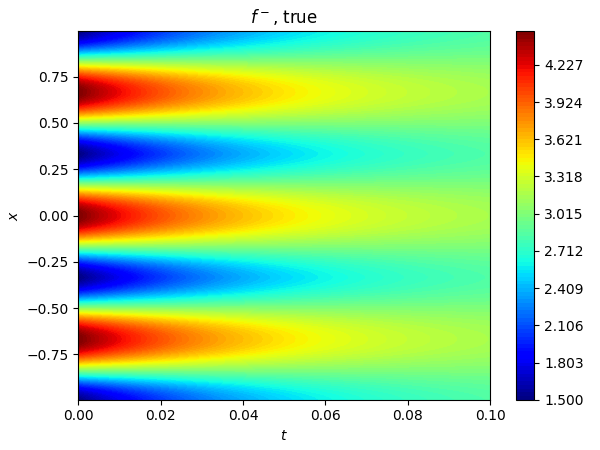}
    \caption{Test 1: Forward problem for the Goldstein-Taylor model with standard PINN in the diffusive regime ($\epsilon=10^{-4}$). Solution of the forward problem by PINN (left) and ground truth (right) of the kinetic densities $f^+$ (top) and $f^-$ (bottom).}
    \label{fig:test1_parabolic_nonap_approximation_forward}
\end{figure}

\begin{figure}[!t]
    \centering
    \includegraphics[width=0.4\textwidth]{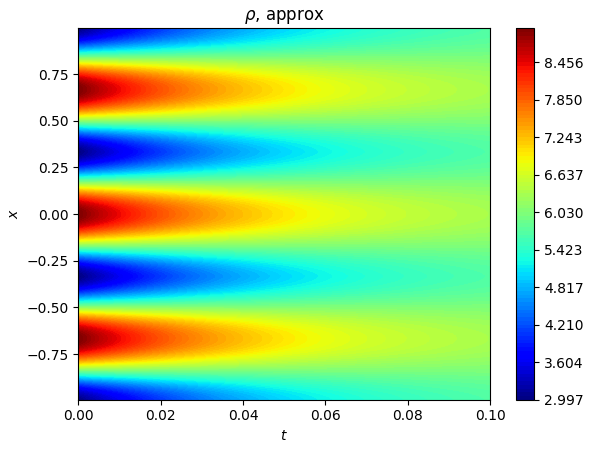}
    \includegraphics[width=0.4\textwidth]{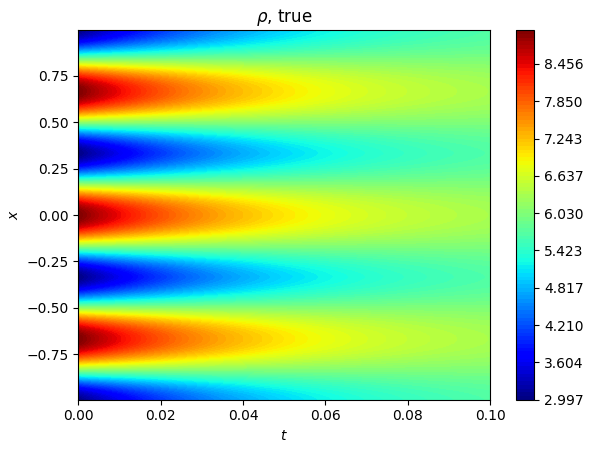}
    \includegraphics[width=0.4\textwidth]{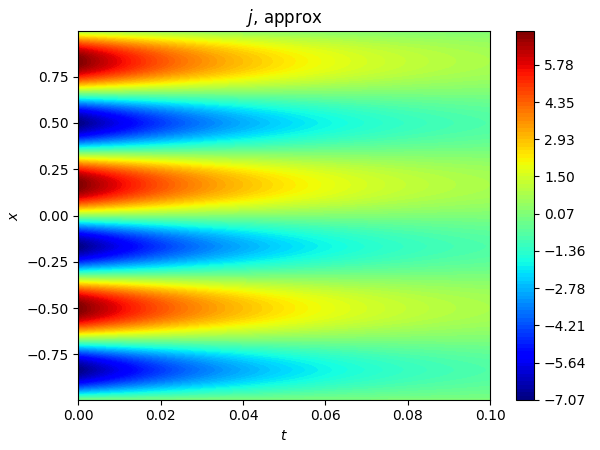} 
    \includegraphics[width=0.4\textwidth] {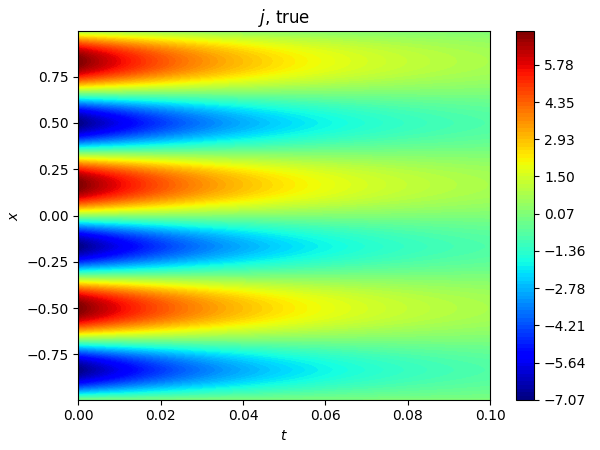}
    \caption{Test 1: Forward problem for the Goldstein-Taylor model with APNN in the diffusive regime ($\epsilon=10^{-4}$). Solution of the forward problem by APNN (left) and ground truth (right) of the density $\rho$ (top) and $j$ (bottom).}
    \label{fig:test1_parabolic_ap_approximation_forward}
\end{figure}
\paragraph{Inverse Problem}

Initially, we consider an inverse problem inferring the scattering coefficient $\sigma$ from the available measurement data using the APNN formulation presented in Appendix \ref{appendix}, with loss function \eqref{eq:pinn-loss_GT} and ${\mathcal L}^\epsilon_r$ term given in \eqref{eq:pinn-MSEr_GT}. For comparison, we also solve the inverse problem applying the standard PINN residual \eqref{eq:pinn-nonap-residue1} in the loss function.
For both APNN and standard PINN formulations, we train the network model on measurements composed of $N_d=24000$ equally spaced samples in the domain $(x, t)\in [-1, 1]\times [0, 0.1]$, from which 20\% (4800) points are randomly selected for validation purpose. For the APNN model we consider measurements only for the density $\rho$, hence assuming to have no information on the flux $j$, whereas for the standard formulation we employ data samples for both the densities $f^+$ and $f^-$ (therefore, in the latter case we assume we have more information on the system \eqref{Goldstein-Taylor}). In addition, $N_r=24000$ residual points are employed with the same data split for validation set. With respect to loss function and training hyperparameters of the APNN given in Table \ref{Table:test12_hyperparameters}, the same setting has been used also for the standard PINN, with the only difference just stated that, when used, the training dataset is given for both variables $f^\pm$ considering equal weights $\omega_d^\pm$.

We show the convergence of the target parameter $\sigma$ in Figure \ref{fig:test1_parabolic_sigma} for both PINN formulations. A very fast convergence can be observed in the APNN, with the initial guess $\sigma_0 = 2$ and a final relative error $\mathcal{O}(10^{-3})$. However, it can be observed that the standard PINN 
failed to recover the correct value of the scattering parameter $\sigma$ (at epoch 4000, early-stopping prevents further training of the PINN). 

\paragraph{Forward Problem}

To further highlight the importance of the AP property, we consider a forward problem for the Goldstein-Taylor model, where scattering coefficient $\sigma=4$ is given and the goal now is to solve the equations on the spatio-temporal domain with corresponding initial conditions. For APNN formulation, $N_b=200$ points are employed to enforce initial conditions of both $\rho$ and $j$, with equation enforced on $N_r=24000$ residual points on the domain $(x, t)\in [-1, 1]\times [0, 0.1]$. The standard PINN formulation  based on the kinetic equations \eqref{Goldstein-Taylor} share the same set with APNN, but initial conditions are given for $f^\pm$.

We plot the solutions obtained with the  standard PINN in Figure \ref{fig:test1_parabolic_nonap_approximation_forward} and with APNN in Figure \ref{fig:test1_parabolic_ap_approximation_forward}. Standard PINN based on the kinetic equations \eqref{Goldstein-Taylor} shows its weakness and converges to a trivial solution on the space-time domain, failing  to approximate the forward solution of both density and flux. On the contrary, the adoption of the APNN ensures the convergence towards the correct diffusive limit, which is also beneficial for the inverse problem we considered before.

\subsection{Test 2: Goldstein--Taylor model with source term}
To examine the performance of the APNN with a more challenging setting  closely related to  epidemic scenarios that we shall discuss later on, we introduce  a source term that creates an oscillatory effect in the density $\rho$ in the Goldstein--Taylor model. The resulting system reads
\begin{equation}
\begin{split}
\frac{\partial \rho}{\partial t} +  \frac{\partial j}{\partial x} &= \kappa \rho\,,	\\ 
\frac{\partial j}{\partial t} + \frac{c^2}{\epsilon^2} \frac{\partial \rho}{\partial x} &= -\frac{\sigma}{\epsilon^2} j \,,
\end{split}
\label{Goldstein-Taylor-source}
\end{equation}
where $\kappa = \kappa(x)$. For this problem, we reformulate the AP-loss function accordingly to the model, simply including the presence of the source term with respect to the formulation discussed in Appendix \ref{appendix}. In the source term, we set $\kappa(x) = \kappa_0 + \kappa_1 \sin(\kappa_2 \pi x)$, with a baseline value $\kappa_0=0$ perturbed by sinusoidal oscillations having amplitude $\kappa_1=3$ and frequency $\kappa_2=4$. We consider again a spatial domain $L=[-1,1]$ and $c=1$. The final goal in this  test is to infer parameters $\kappa_0$, $\kappa_1$ and $\kappa_2$ and evaluate the spatio-temporal reconstruction given by the APNN with a partially observed system, having only information of $\rho$, considering a scattering coefficient $\sigma = 1$ and the following initial conditions:
$$\rho(x,0) = 1 + 0.5\,e^{-10 x^2}, \qquad j(x,0) = 10 x \,e^{-10 x^2}.$$
\begin{table}[!b]
    \centering
    \begin{tabular}{@{}ccccc@{}}\toprule 
    Parameter & Ground Truth & Initial Guess & Estimation & Relative Error \\
    \midrule 
    $\kappa_0$ & 0 & 0.5 & 0.0011 & N/A \\
    $\kappa_1$ & 3 & 2 & 2.9263 & $2.46\times 10^{-2}$ \\
    $\kappa_2$ & 4 & 3 & 4.0003 & $7.50\times 10^{-5}$ \\
    \bottomrule
    \end{tabular}
    \caption{Test 2 (a): Goldstein-Taylor model with source in diffusive regime ($\epsilon=10^{-5}$) with density data only. Inference results for the source term coefficients using the APNN.}
    \label{Table:test2_GaussianIC-parabolic}
\end{table}
\begin{figure}[!t]
    \centering
    \includegraphics[width=0.32\textwidth]{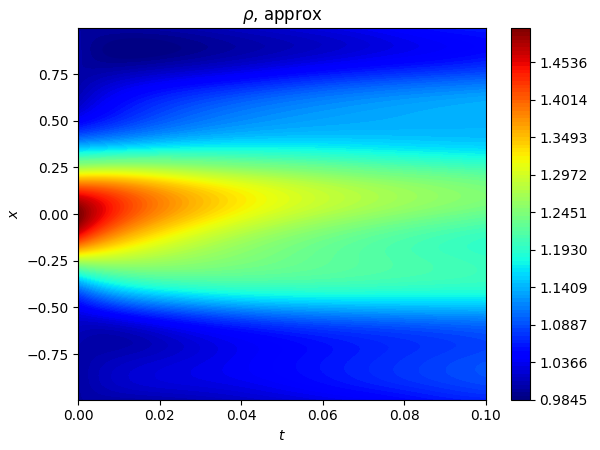}
    \includegraphics[width=0.32\textwidth]{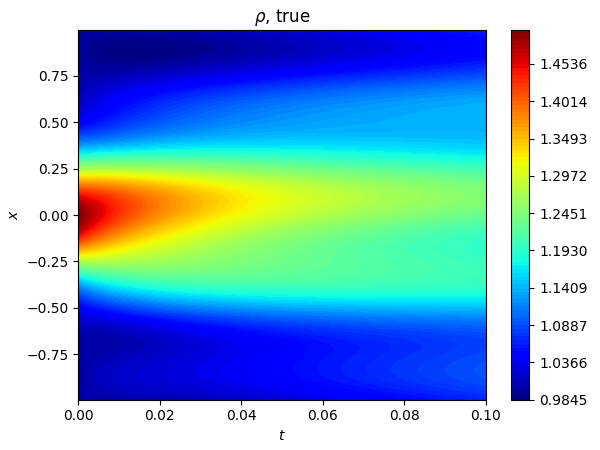}
    \includegraphics[width=0.32\textwidth]{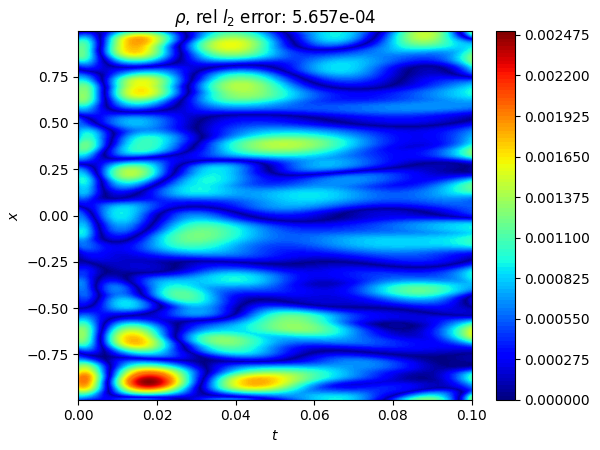}
    \includegraphics[width=0.32\textwidth]{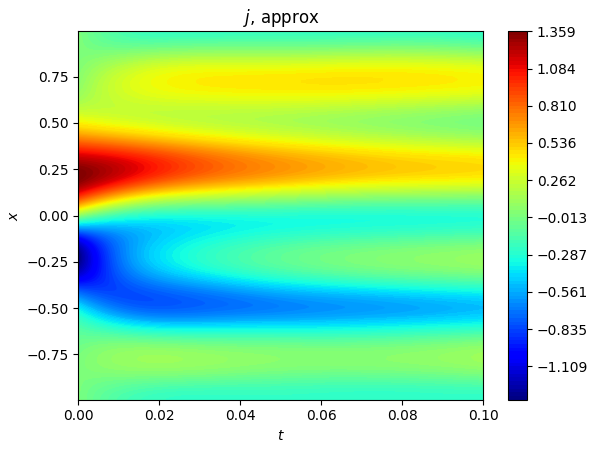}
    \includegraphics[width=0.32\textwidth]{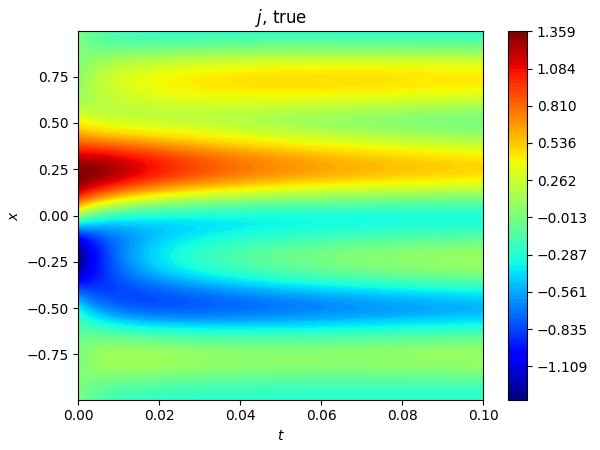}
    \includegraphics[width=0.32\textwidth]{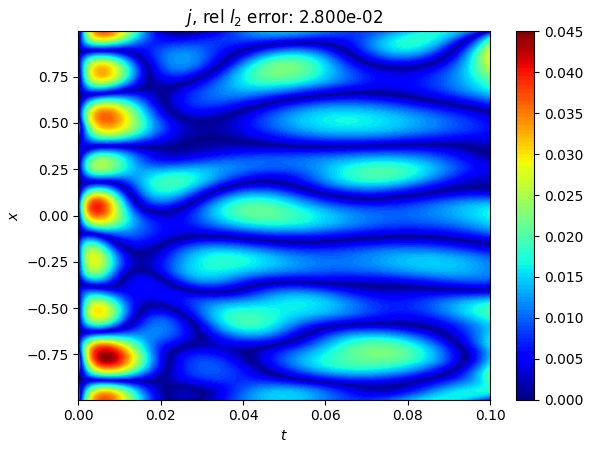}
    \caption{Test 2 (a): Goldstein-Taylor model with source in diffusive regime ($\epsilon=10^{-5}$) with density data only. Approximated forward solution (left column), ground truth (middle column) and relative $L^2$ error (right column) of density $\rho$ (first row) and flux $j$ (second row) obtained with the APNN.}
    \label{fig:test2_GaussianIC-parabolic}
\end{figure}
\subsubsection*{Test 2 (a): Diffusive regime with density data only}
We initially consider a diffusive, parabolic regime defined by $\epsilon=10^{-5}$, with $t_{end}=0.1$. We employ $N_d = 12000$ for $\rho$, not considering any dataset for $j$, while still imposing initial and boundary conditions for both variables. For the residual term, we use $N_r = 12000$ points on the domain $(x, t)\in [-1, 1]\times [0, 0.1]$. We use 20\% of $N_d$ and $N_r$ for validation purposes and the rest for the training. 
Results of the parameters inference are shown in Table \ref{Table:test2_GaussianIC-parabolic}, where initial guesses of target variables are listed, even though we observed that the neural network is not very  sensitive to the choice of these values. From these results we can observe that, in general, the most difficult coefficient to calibrate with the NN is the amplitude of the perturbation of the source term, $\kappa_1$.

The APNN forward approximations of $\rho$ and $j$ are presented in Figure \ref{fig:test2_GaussianIC-parabolic}, where we can observe that forward solutions well capture the correct dynamics of $\rho$ and accurately recover  $j$  without any measurement on the latter. Nonetheless, we acknowledge that when concerning diffusive regimes as in Eq. \eqref{eq.heat}, the problem results fully described by the sole density $\rho$, and the absence of information on $j$ does not lead to an actual lack of data knowledge.
\begin{table}[!b]
    \centering
    \begin{tabular}{@{}ccccc@{}}\toprule 
    Parameter & Ground Truth & Initial Guess & Estimation & Relative Error \\
    \midrule 
    $\kappa_0$ & 0 & 0.5 & $4.37\times 10^{-5}$ & N/A \\
    $\kappa_1$ & 3 & 2 & 3.0005 & $1.67\times 10^{-4}$ \\
    $\kappa_2$ & 4 & 3 & 4.0002 & $5.00\times 10^{-5}$ \\
    \bottomrule
    \end{tabular}
    \caption{Test 2 (b): Goldstein-Taylor model with source in hyperbolic regime ($\epsilon=1$) with density data only. Inference results for the source term coefficients using the APNN.}
    \label{Table:test2_GaussianIC-hyperbolic}
\end{table}
\begin{figure}[!t]
    \centering
    \includegraphics[width=0.32\textwidth]{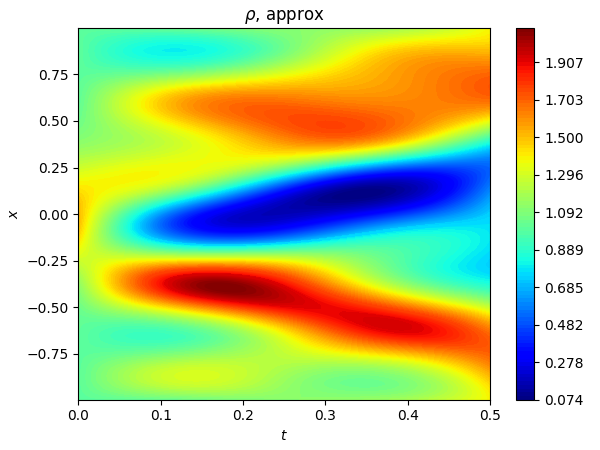}
    \includegraphics[width=0.32\textwidth]{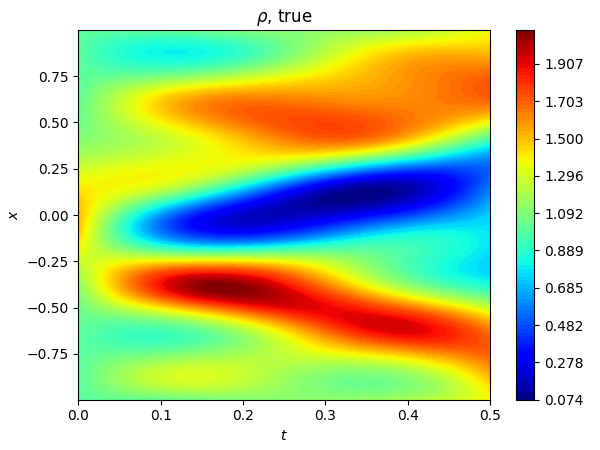}
    \includegraphics[width=0.32\textwidth]{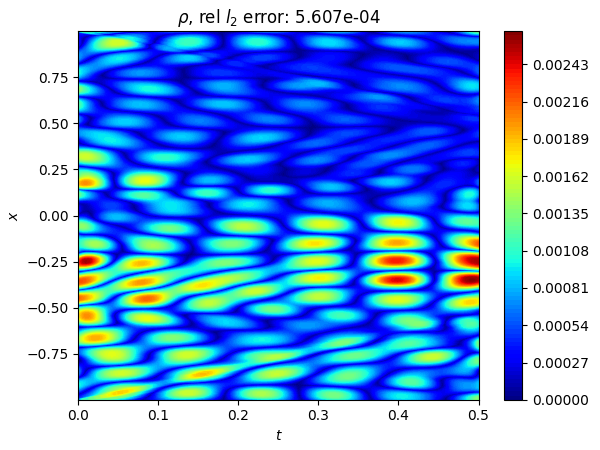}
    \includegraphics[width=0.32\textwidth]{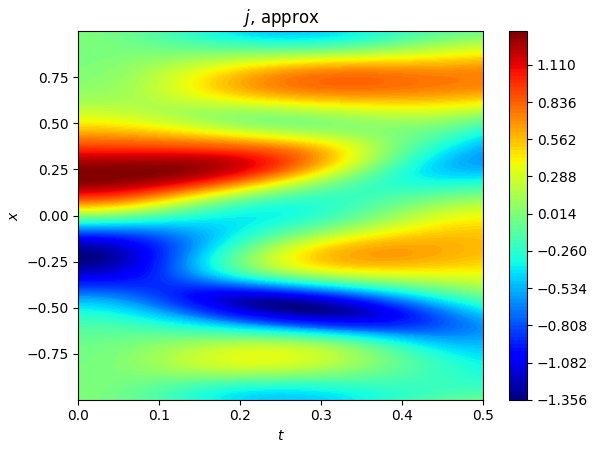}
    \includegraphics[width=0.32\textwidth]{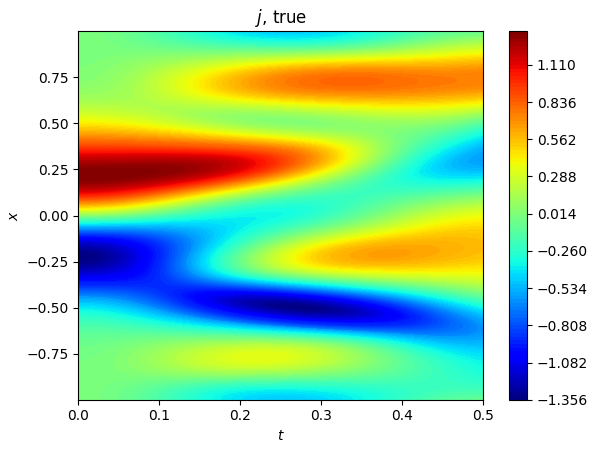}
    \includegraphics[width=0.32\textwidth]{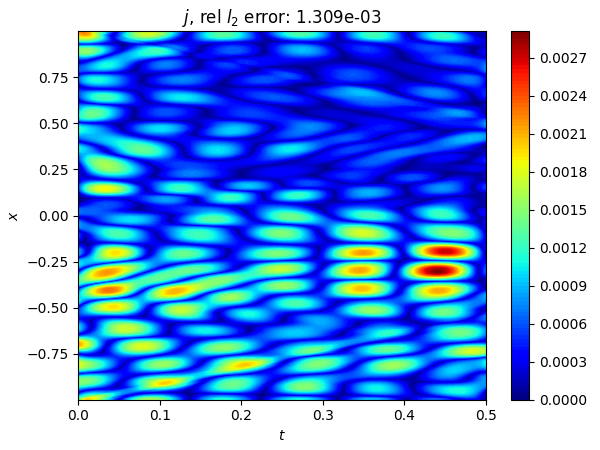}
    \caption{Test 2 (b): Goldstein-Taylor model with source in hyperbolic regime ($\epsilon=1$) with density data only. Approximated forward solution  (left column), ground truth (middle column) and relative $L^2$ error (right column) of density $\rho$ (first row) and flux $j$ (second row) obtained with the APNN.}
    \label{fig:test2_GaussianIC-hyperbolic}
\end{figure}
\subsubsection*{Test 2 (b): Hyperbolic regime with density data only}
In the second case, we consider a hyperbolic regime with $\epsilon=1$ and $t_{end}=0.5$. We employ $N_d = 16800$ for $\rho$, not considering again any dataset for $j$, and fix $N_r = 16800$ on the domain $(x, t)\in [-1,1]\times [0, 0.5]$, with 20\% of each dataset for validation. Coefficients inferred by the APNN are listed in Table  \ref{Table:test2_GaussianIC-hyperbolic}, while forward solutions are shown in Figure  \ref{fig:test2_GaussianIC-hyperbolic}. 
Similar to the diffusive regime, the APNN correctly infer all the unknown parameters and is capable of  approximating the solution of densities $\rho$ and $j$ well, but in this case in a much more demanding problem. Indeed, even though in hyperbolic regimes the problem is not completely defined by the sole density of the system, being the dataset really incomplete without any information on the flux $j$, the APNN is still capable of approximating the correct solution of the whole dynamics.
\subsection{Test 3: SIR transport model with constant epidemic parameters}
In the following, we evaluate the performance of the APNN with respect to the dynamics governed by the SIR multiscale transport model \eqref{eq.SIR_kinetic}.
We first design a numerical test with an initial condition that simulates the presence of two epidemic hot-spots, aligned in the spatial domain $L=[0,20]$, presenting a different number of infected individuals, distributed following a Gaussian function,
$$ I(x,0)= \alpha_1\,e^{-(x-x_1)^2} + \alpha_2\,e^{-(x-x_2)^2}, $$
where $x_1=5$ and $x_2=15$ are the coordinates of the hot-spots, while 
$\alpha_1=0.01$ and $\alpha_2=0.0001$ define the different initial epidemic concentration in the two cities, hence with a deeply higher density of infected individuals in the first city.
Assuming that there are no immune individuals at $t=0$ and that the total population is uniformly distributed in the domain, we have
$$ S(x,0)=1 - I(x,0), \qquad R(x,0)=0. $$
We impose initial fluxes in equilibrium, following \eqref{eq:fick}, and periodic boundary conditions to allow both directions of connection for the two cities.
We initially consider a simple setting defined by constant epidemic parameters in space and time, with $\beta=12$ and $\gamma = 6$, which lead to study an infectious disease characterized by an initial reproduction number $R_t(0) = 2$. 

The APNN is used to infer both the epidemic parameters as well as approximate the solutions for a parabolic and a hyperbolic scenario. To mimic the availability of data close to reality, we use a sparse dataset for the training process, sampling the spatio-temporal points from the available dataset with probability proportional to the magnitude of $I$. We consider, indeed, that in real-world epidemic scenarios data on the evolution of the infectious disease are only available in the regions in which the virus has already started to spread. Specifically, the probability of each spatio-temporal location $(x, t)$ chosen for the training dataset is given by
\begin{equation}
\label{eq:importance_sampling}
    p(x, t) = \frac{I(x, t)}{\int_{\Omega}I(x, t)}.
\end{equation}
\begin{figure}[!t]
    \centering   
    \includegraphics[width=0.32\textwidth]{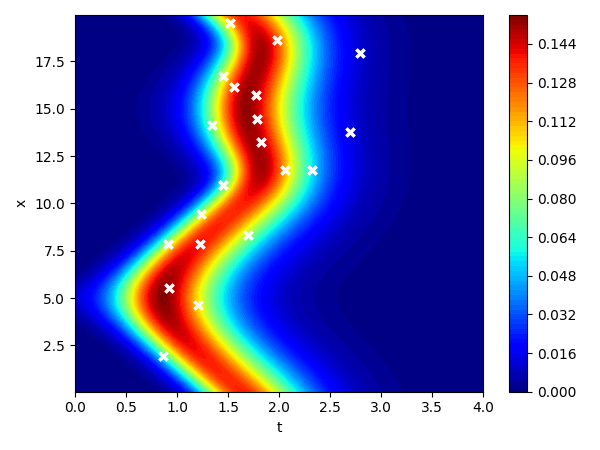}
    \includegraphics[width=0.32\textwidth]{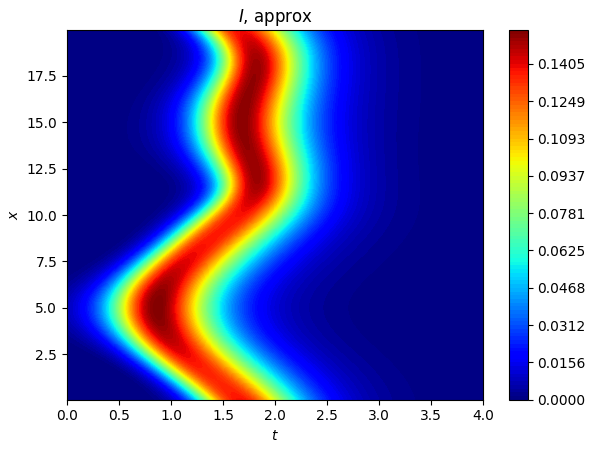}
    \includegraphics[width=0.32\textwidth]{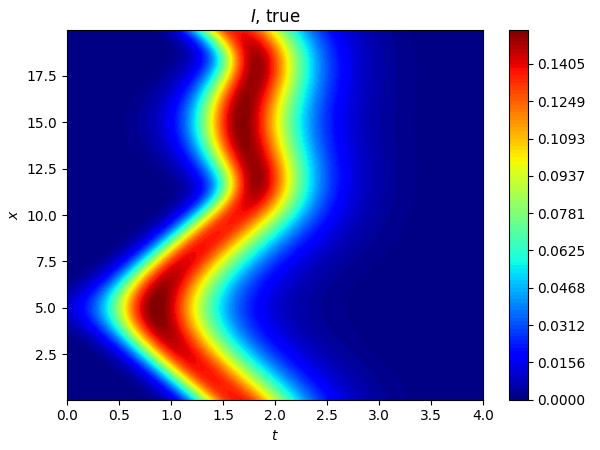}
    \caption{Test 3.1 (a): SIR transport model with constant epidemic parameters and partially observed dynamics in diffusive regime ($\lambda_{S,I,R}^2=10^3$, $\tau_{S,I,R}=10^{-3}$). Selected sparse samples ($N_d=20$) marked with white crosses (left column), approximation obtained in the inverse problem (middle column), and ground truth (right column) of the densities of infected $I$.  }
    \label{fig:test3_parabolic_inverse}
\end{figure}

\begin{figure}[!t]
    \centering   
    \includegraphics[width=0.32\textwidth]{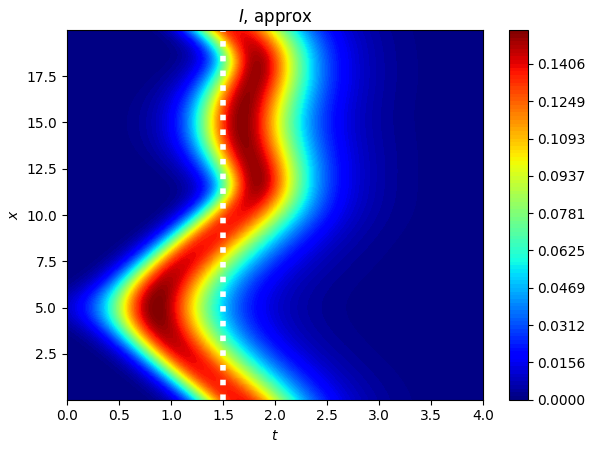}
    \includegraphics[width=0.32\textwidth]{figures_final/test2/scattered_data_corrected/test2_parabolic_scattered_I_true}
    \includegraphics[width=0.32\textwidth]{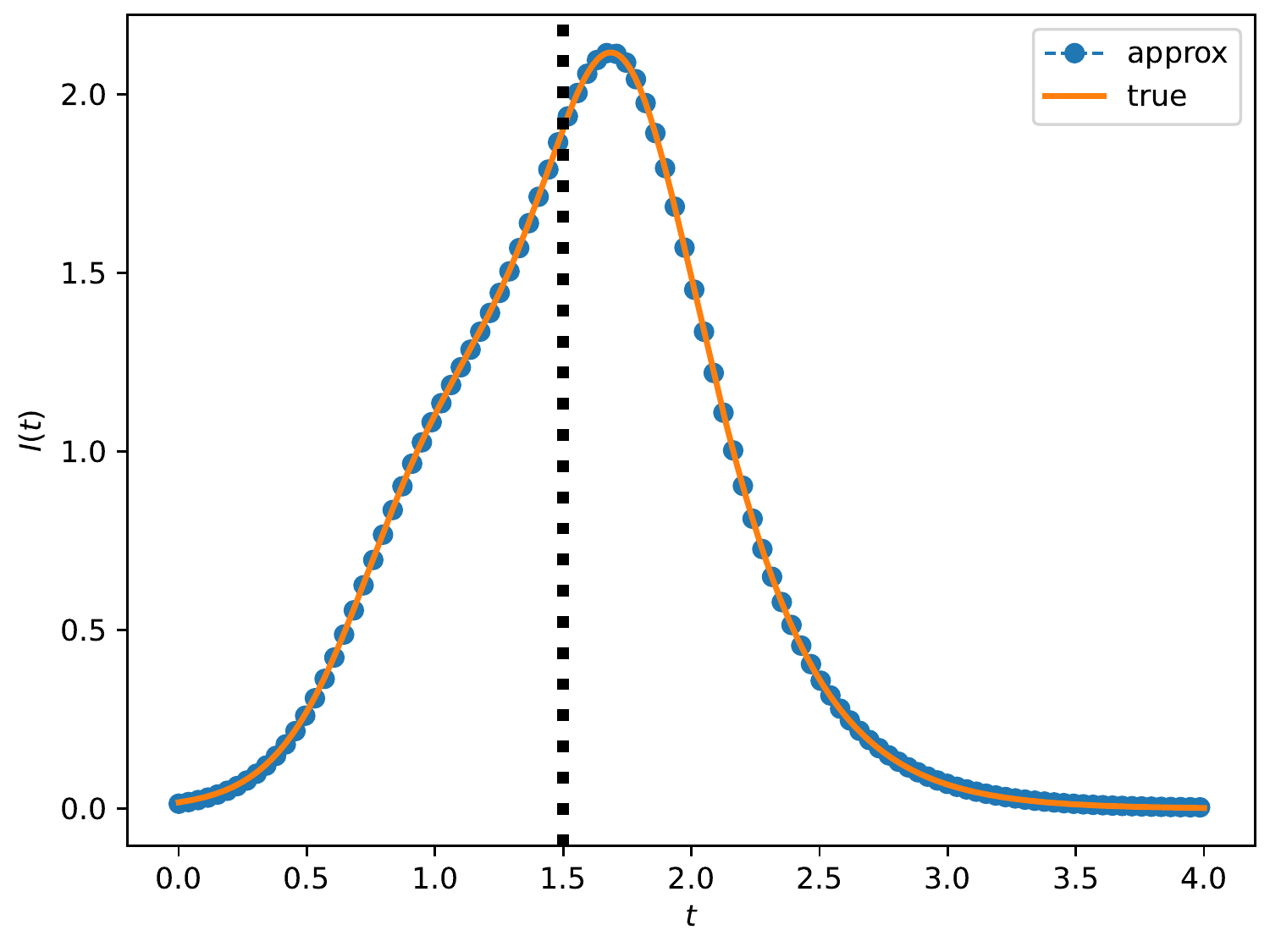}
    \includegraphics[width=0.32\textwidth]{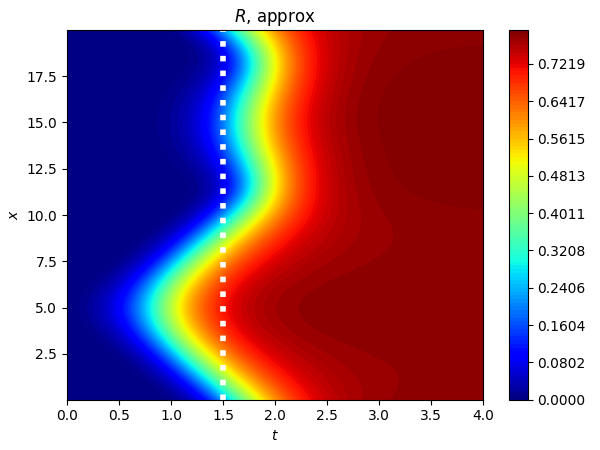}
    \includegraphics[width=0.32\textwidth]{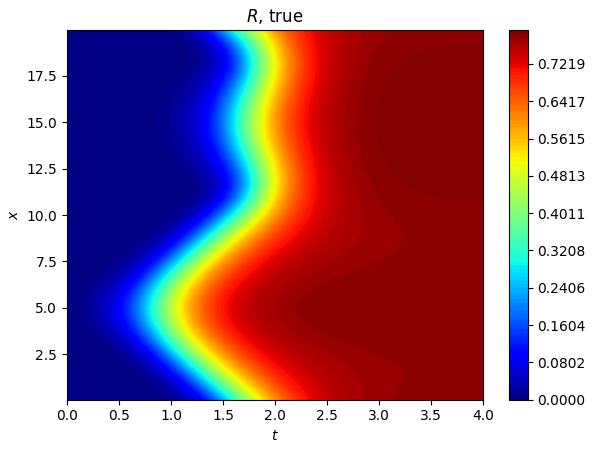}
    \includegraphics[width=0.32\textwidth]{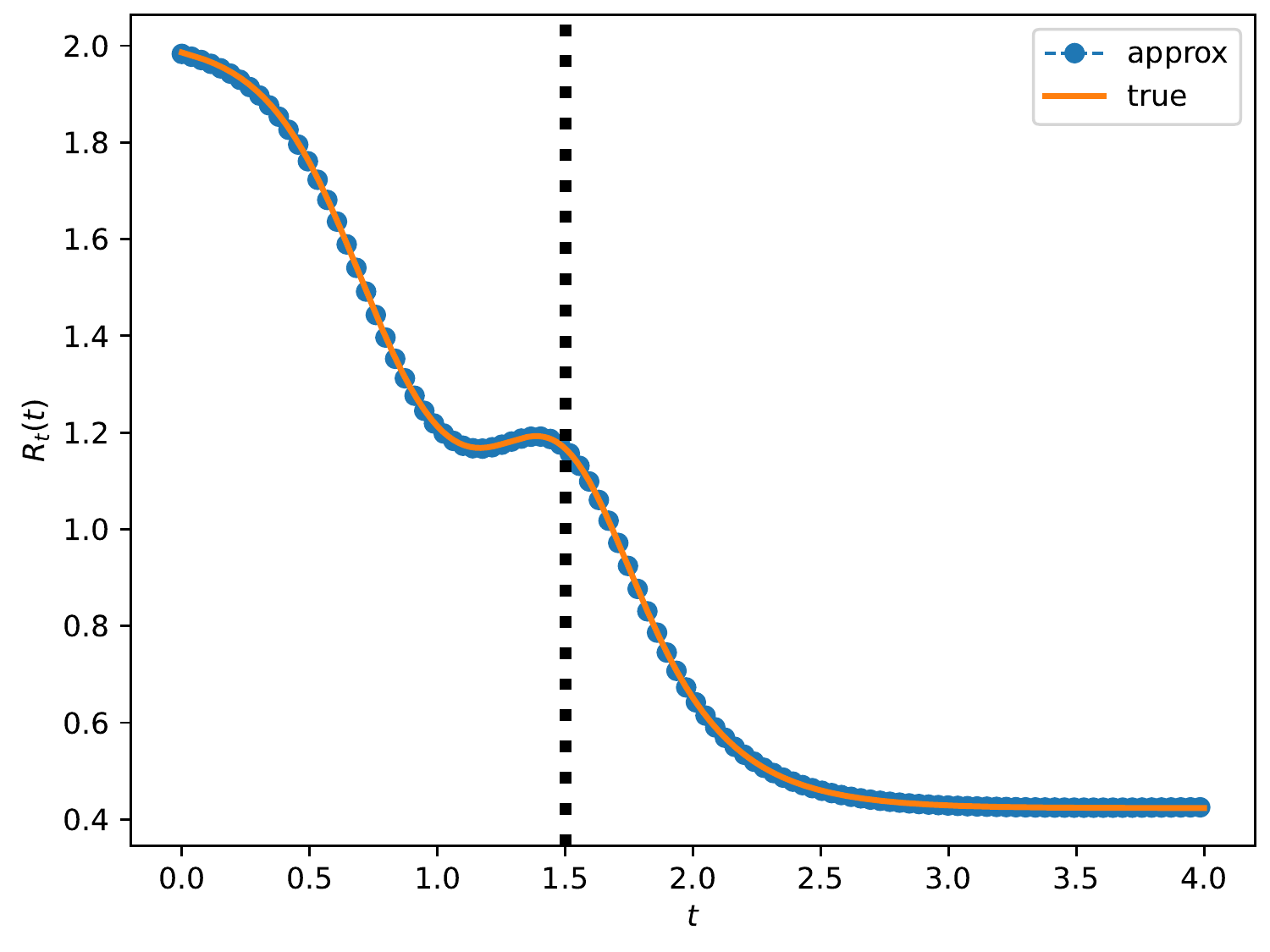}
    \caption{Test 3.2 (a): SIR transport model with constant epidemic parameters and partially observed dynamics in diffusive regime ($\lambda_{S,I,R}^2=10^3$, $\tau_{S,I,R}=10^{-3}$). Approximation and forecast with measurements on a short time $t \in [0,1.5]$ denoted by the dashed line (left column), and ground truth (middle column) of infected $I$ (first row) and removed $R$ (second row). Temporal evolution of the cumulative density of infected individuals $I$ in the whole domain (first row, right) and of the reproduction number $R_t$ (second row, right) obtained with the APNN, trained  based on a short time period $t \in [0,1.5]$ (marked by the dotted line).  }
    \label{fig:test3_parabolic_forward}
\end{figure}
\subsubsection*{Test 3 (a): Partially observed dynamics in diffusive regime}
In the first case, a parabolic configuration of speeds and relaxation parameters is considered, with $\lambda_{S,I,R}^2=10^3$ and $\tau_{S,I,R}=10^{-3}$. We examine the performance of the APNN in the two following different problems.
\begin{itemize}
\item {\bf Test 3.1 (a): Parameter inference test.} We consider a sparse dataset where only $N_d=20$ measurements are selected from the entire space-time domain  $(x,t) \in [0,20] \times [0,4]$, according to the density of $I$, as described in \eqref{eq:importance_sampling}, as shown in Figure \ref{fig:test3_parabolic_inverse} (left).
\item {\bf Test 3.2 (a): Forecasting test.} As a second problem, we intend to investigate the forecasting capability of the APNN. In contrast with sampling measurements available across the entire spatio-temporal domain in the parameter inference test, we generate a training dataset of size $N_d=5300$ on a shorter time domain $t \in [0,1.5]$ 
and we assess the correctness of APNN approximations in $t \in [0,1.5]$ and forecasting performance in $t \in [1.5,4]$.
\end{itemize} 
In both cases, equations residual are enforced on $N_r=40000$ residual points on the spatio-temporal domain and 20\% of each dataset is used for validation. In addition, we assume initial conditions for $S,I,R$ are unknown in both problems, thus requiring an even more demanding performance to the APNN. 

\begin{table}[!b]
    \centering
    \begin{tabular}{@{}ccccc@{}}\toprule 
    Parameter & Ground Truth & Initial Guess & Estimation & Relative Error \\
    \midrule 
    $\beta$ & 12 & 8 & 11.9428 & $4.76\times 10^{-3}$ \\
    $\gamma$ & 6 & 3 & 5.9772 & $3.80\times 10^{-3}$ \\
    \bottomrule
    \end{tabular}
    \caption{Test 3.1 (a): SIR transport model with constant epidemic parameters and partially observed dynamics in diffusive regime ($\lambda_{S,I,R}^2=10^3$, $\tau_{S,I,R}=10^{-3}$). Inferred results for transmission rate $\beta$ and recovery rate $\gamma$ from a sparse measurement dataset of $N_d=20$ samples, and the relative error with respect to the ground truth values.}
    \label{Table:test3_parabolic_inference}
\end{table}

Results of the parameter inference  task based on the sparse measurement dataset are reported in Table \ref{Table:test3_parabolic_inference}, where an excellent estimation of both $\beta$ and $\gamma$ can be observed with respect to the ground truth. 
Figure \ref{fig:test3_parabolic_inverse} shows that the reconstructed forward approximations for the density of the epidemic compartment $I$ have an excellent agreement with the  true solution in the entire  domain. 
Also in the forecasting test,  the approximated and predicted dynamics (based on the measurements from the time period $t \in [0, 1.5]$) perfectly match the ground truth in the entire domain $t \in [0,4]$, as shown in Figure \ref{fig:test3_parabolic_forward}, even if in this demanding setting initial conditions of densities are assumed to be unknown. 
These results further highlight the capability of APNN to forecast the spread of an infectious disease in diffusive regimes thanks to the physical knowledge of the PDE system embedded in the NN together with the preservation of the AP property.
In the same Figure, we present also the temporal evolution of the cumulative density of infected individuals $I$ in the whole domain as well as the effective reproduction number $R_t$ predicted by the APNN. The excellent agreement between predictions ($t>1.5$) and the ground truth out of the training domain further assess the forecasting capability of APNNs.
\begin{table}[!b]
    \centering
    \begin{tabular}{@{}ccccc@{}}\toprule 
    Parameter & Ground Truth & Initial Guess & Estimation & Relative Error \\
    \midrule 
    $\beta$ & 12 & 8 & 12.0126 & $1.05\times 10^{-3}$ \\
    $\gamma$ & 6 & 3 & 6.0447 & $7.45\times 10^{-3}$ \\
    \bottomrule
    \end{tabular}
    \caption{Test 3.1 (b): SIR transport model with constant epidemic parameters and partially observed dynamics in hyperbolic regime ($\lambda_{S,I,R}^2=1$, $\tau_{S,I,R}=1$). Inferred results for transmission rate $\beta$ and recovery rate $\gamma$ from a sparse measurement dataset of $N_d=20$ samples, and relative error with respect to the ground truth values.}
    \label{Table:test3_hyperbolic_inference}
\end{table}
\begin{figure}[!t]
    \centering   
    \includegraphics[width=0.32\textwidth]{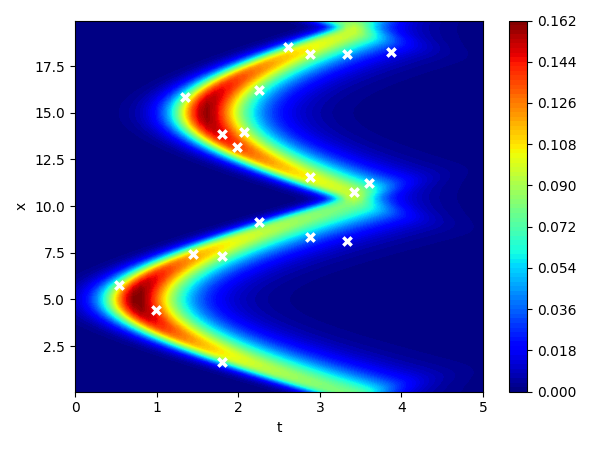}
    \includegraphics[width=0.32\textwidth]{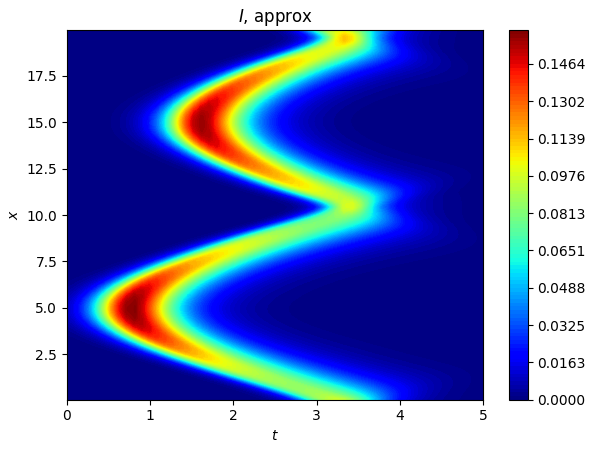}
    \includegraphics[width=0.32\textwidth]{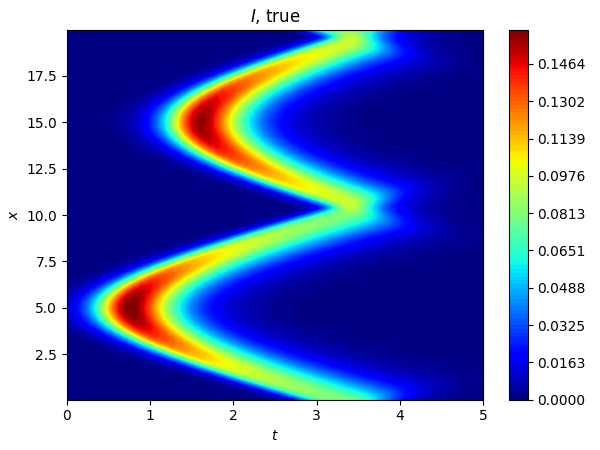}
    \caption{Test 3.1 (b): SIR transport model with constant epidemic parameters and partially observed dynamics in hyperbolic regime ($\lambda_{S,I,R}^2=1$, $\tau_{S,I,R}=1$). Selected sparse samples ($N_d=20$) marked with white crosses (left column), approximation obtained in the inverse problem (middle column), and ground truth (right column) of the densities of infected $I$. }
    \label{fig:test3_hyperbolic_inverse}
\end{figure}

\begin{figure}[!t]
    \centering   
    \includegraphics[width=0.32\textwidth]{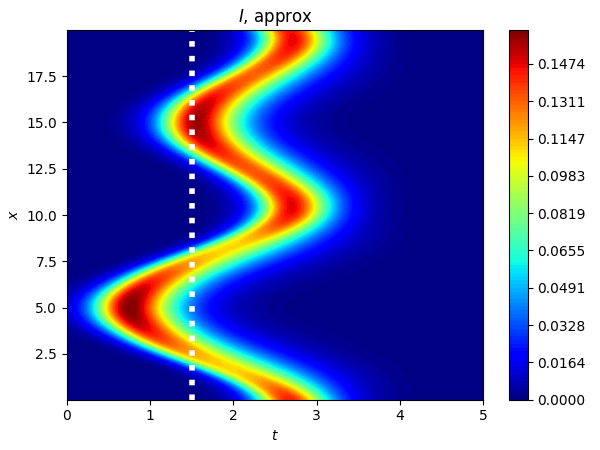}
    \includegraphics[width=0.32\textwidth]{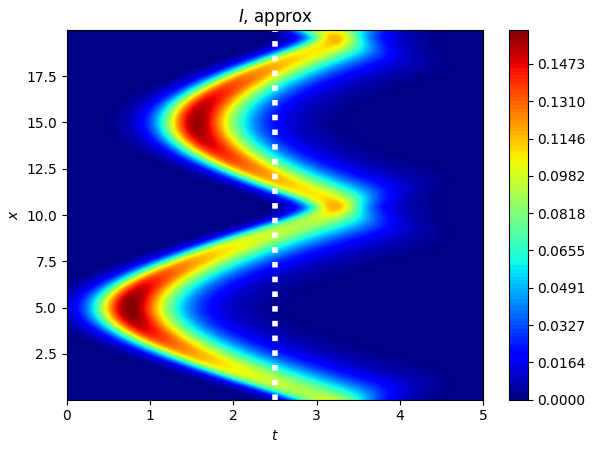}
    \includegraphics[width=0.32\textwidth]{figures_final/test2/scattered_data_corrected/test2_hyperbolic_scattered_I_true}
    \includegraphics[width=0.32\textwidth]{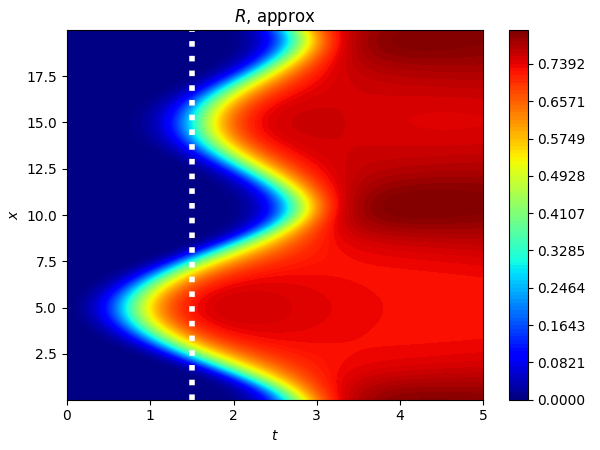}
    \includegraphics[width=0.32\textwidth]{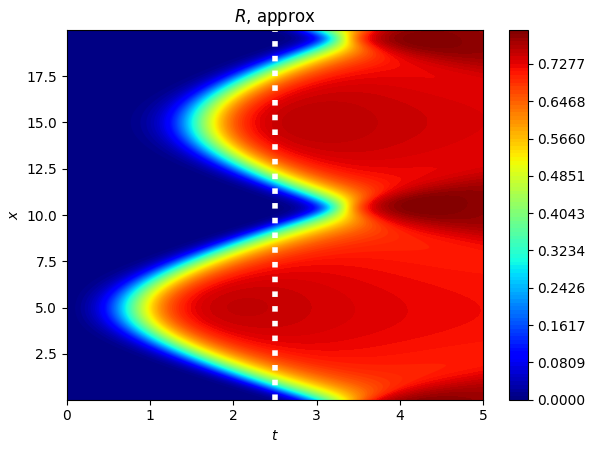}
    \includegraphics[width=0.32\textwidth]{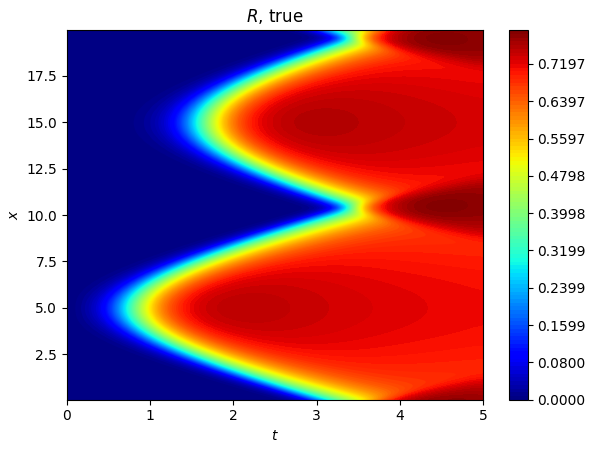}
    \caption{Test 3.2 (b): SIR transport model with constant epidemic parameters and partially observed dynamics in hyperbolic regime ($\lambda_{S,I,R}^2=1$, $\tau_{S,I,R}=1$). Approximation with measurements from a shorter time period $t \in [0,1.5]$ (first column) or $t \in [0,2.5]$ (middle column), stopped at the dashed line, and ground truth (last column), of the densities of infected $I$ (first row) and removed $R$ (second row).     }
    \label{fig:test3_hyperbolic_forward}
\end{figure}

\begin{figure}[!t]
    \centering
    \includegraphics[width=0.45\textwidth]{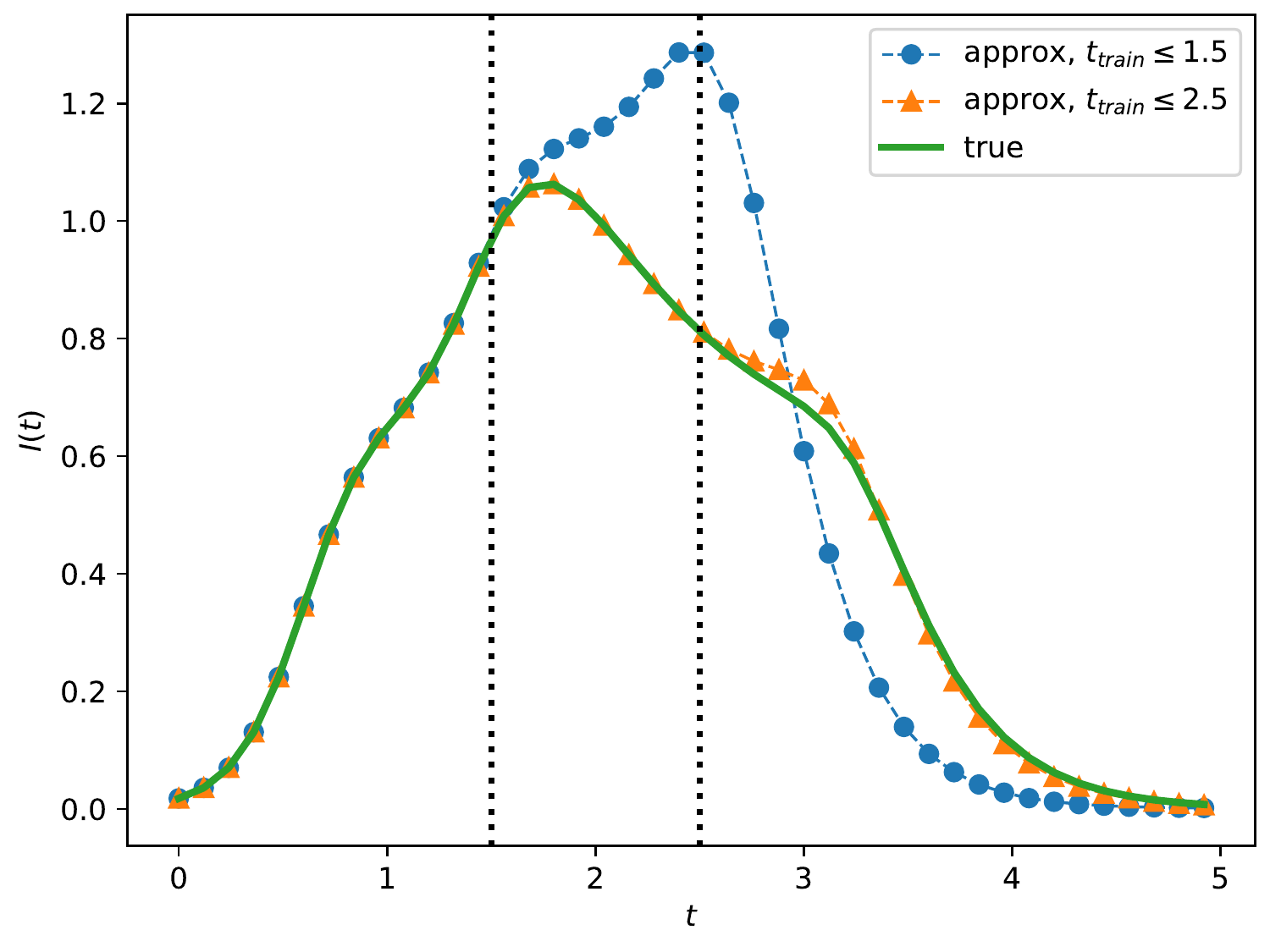}\hskip .5cm
    \includegraphics[width=0.45\textwidth]{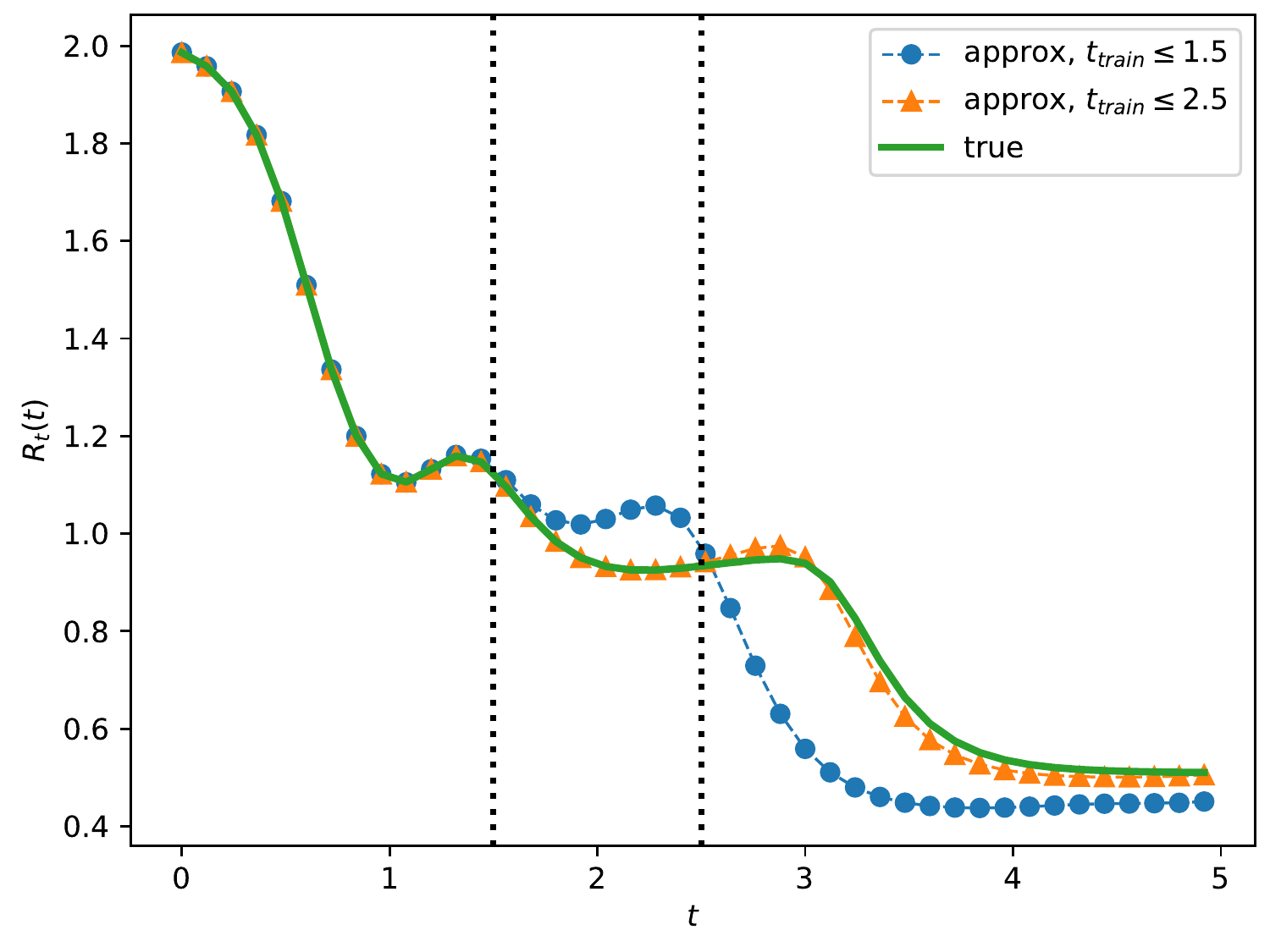}
    \caption{
    Test 3.2 (b), SIR transport model with constant epidemic parameters and partially observed dynamics in hyperbolic regime ($\lambda_{S,I,R}^2=1$, $\tau_{S,I,R}=1$). Temporal evolution of the cumulative density of infected individuals $I$ in the whole domain (left) and of the reproduction number $R_t$ (right) obtained with the APNN using measurements from a shorter period of $t \in [0,1.5]$ or $t \in [0,2.5]$ (stopped at the dotted lines) compared with ground truth.    }
    \label{fig:test3_hyperbolic_Rt}
\end{figure}
\subsubsection*{Test 3 (b): Partially observed dynamics in hyperbolic regime}
In the second case, we consider a hyperbolic regime with $\lambda_{S,I,R}=1$ and $\tau_{S,I,R}=1$.  
As previously done, we consider two different contexts.
\begin{itemize}
\item  {\bf Test 3.1 (b): Parameter inference test.} We first consider a sparse measurement setting, where $N_d = 20$ measurements the spatio-temporal domain $(x, t)\in [0, 20]\times [0, 5]$ are available. The chosen samples are shown in Figure \ref{fig:test3_hyperbolic_inverse} (left) and have been selected again according to the density of $I$, as described in \eqref{eq:importance_sampling}.
\item  {\bf Test 3.2 (b): Forecasting test.} Secondly, we consider a forecasting task, training the APNN with the measurements generated from a limited time domain $t \in [0,t_{train}]$. In this example, we chose $t_{train}=1.5$ and $t_{train}=2.5$ with $N_d = 5000$ and $N_d = 8500$ measurements of densities employed respectively, and then evaluate the network performance over the time domain $t \in [0,5]$. 
\end{itemize}
In both scenarios, $N_r=23600$ residual points are employed on the spatio-temporal domain to enforce the underlying equations, still assuming that initial conditions of densities $S,I,R$ are unknown, as in the previous test case. 

Parameters $\beta$ and $\gamma$ estimated by the APNN from the sparse measurements are presented in Table \ref{Table:test3_hyperbolic_inference}, where we observe again a very good agreement with respect to true values.  At the same time, the APNN is capable of reconstructing the correct dynamics of the phenomenon of interest in the whole domain besides the sparsity and incompleteness of data, as shown in Figure~\ref{fig:test3_hyperbolic_inverse}.
On the other hand, we show results obtained when training the APNN with measurements taken from a shorter time period in Figures \ref{fig:test3_hyperbolic_forward} and \ref{fig:test3_hyperbolic_Rt}. 
In Figure \ref{fig:test3_hyperbolic_Rt}, we plot the temporal evolution of the cumulative density of infected individuals $I$ in the whole domain as well as the effective reproduction number $R_t$ predicted by the APNN when the measurement data restrict to the shorter time periods $t \in [0,t_{train}]$. When $t_{train}=1.5$, APNN predictions deviate from the ground truth almost immediately after the training period. In contrast, when the measurement data is extended to $t_{train}=2.5$, the APNN produces good reconstructions and predictions in the forecasting region ($t>2.5$) with respect to the ground truth. Similar observations can be  made for the approximations of densities $I$ and $R$ presented in Figure \ref{fig:test3_hyperbolic_forward}. 
This behavior of the APNN is observed because, when considering a dataset only for $t \in [0,1.5]$, we do not cover enough information of the major dynamics. We remark indeed that no data on fluxes is given to the APNN, which, in a hyperbolic regime, means to deal with a consistent lack of data knowledge. The predictions obtained, in fact, show that the APNN tends to smooth out the actual epidemic propagation pattern, not describing the correct transport/hyperbolic mechanism in the regions connecting the two urban areas. This appears clear when looking at Figure \ref{fig:test3_hyperbolic_forward} (first column) and Figure \ref{fig:test3_hyperbolic_Rt}, and observing that the dynamics predicted by the APNN tend to spread the virus faster in a more diffusive way, giving rise to a fake epidemic hot-spot around $t = 2.5$. 
\begin{figure}[!t]
\centering
\begin{minipage}{6in}
  \centering
  \raisebox{-0.5\height}{\includegraphics[height=4.5 cm]{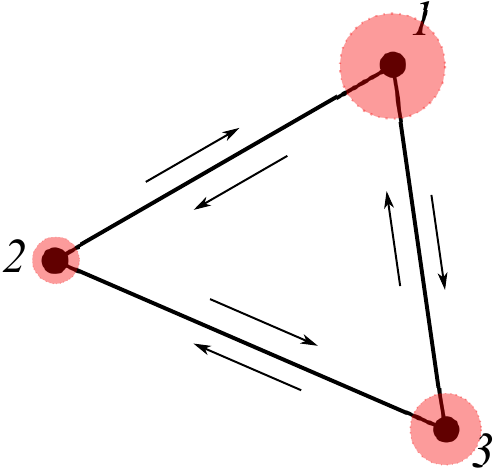}}
  \hspace*{.2in}
  \raisebox{-0.5\height}{\includegraphics[height=7 cm]{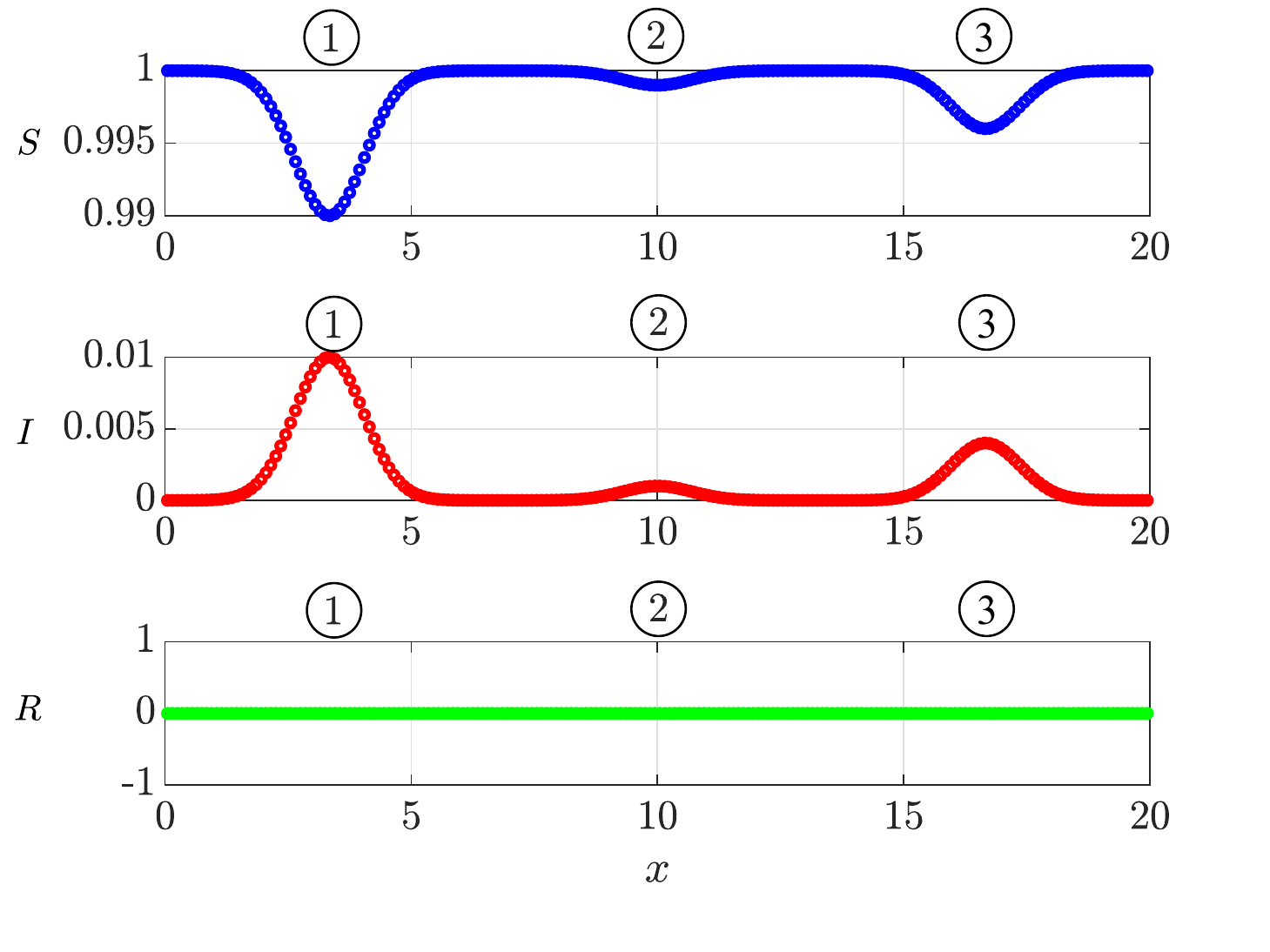}}
  \vspace*{-.2in}
\end{minipage}
\caption{Test 4: SIR transport model with spatially variable transmission rate. Left: schematic representation of the spatial setting considered, with 3 initial hot-spots presenting different initial concentrations of infectious people, proportional to the light red circles. 
Individuals move from one location to another following the two opposite directions defined in the one-dimensional space with periodic boundary conditions. Due to the heterogeneous environment 3 additional hot spots will form along the main connection lines.  Right: initial conditions for susceptible (top), infectious (middle) and removed (bottom).}
\label{fig.scheme_test2}
\end{figure}

\subsection{Test 4: SIR transport model with heterogeneous environment}
\label{section_test4}
Next, we consider a much more challenging scenario, taking into account a spatially varied transmission rate that follows a hypothetical heterogeneous environment. An initial condition of the SIR multiscale transport model is designed to simulate the presence of 3 epidemic hot-spots aligned in the spatial domain $L=[0,20]$, each one having a different initial density of infected individuals, distributed in space following again a Gaussian:
$$ I(x,0)= \alpha_1\,e^{-(x-x_1)^2} + \alpha_2\,e^{-(x-x_2)^2} + \alpha_3\,e^{-(x-x_3)^2}.$$
Here $x_1=10/3$, $x_2=10$, $x_3=50/3$ are the coordinates of the epidemic centers and 
$\alpha_1=0.01$, $\alpha_2=0.001$, $\alpha_3=0.004$ define the different initial epidemic concentration in each spot.
Assuming again that there are no immune individuals at $t=0$ and that the total population is uniformly distributed in the spatial domain, we set $S(x,0)=1 - I(x,0)$ and $R(x,0)=0$.
As previously, we impose initial fluxes in equilibrium, following \eqref{eq:fick}, and periodic boundary conditions, to allow a connection also between hot-spots 1 and 3, so that the domain connecting the positions of these regions form the closed shape presented in Figure \ref{fig.scheme_test2} (left). In the same figure  (right), initial conditions of the 3 epidemic compartments are shown.
We set the following spatially varied transmission rate^^>\cite{Bert,Wang2020}:
\begin{equation*}
\beta(x) = \beta_0 + \beta_1\sin\left(\zeta\pi x\right),
\end{equation*}
with $\beta_0 = 9$ and perturbing this baseline value with oscillations of amplitude $\beta_1 = 2.5$ and frequency $\zeta = 0.55$. The recovery rate is set to be $\gamma = 8$.
This choice of  parameters  simulates an infectious disease characterized by an initial reproduction number $R_t(0) \approx 1.05$.
With the APNN, the goal is to infer $\beta_0, \beta_1$ and $\zeta$ as well as approximate the dynamics of densities based on the partially available measurements and the forecasting performance.

\begin{table}[!b]
    \centering
    \begin{tabular}{@{}ccccc@{}}\toprule 
    Parameter & Ground Truth & Initial Guess & Estimation & Relative Error \\
    \midrule 
    $\beta^0$ & 9 & 5 & 9.0170 & $1.89\times 10^{-3}$ \\
    $\beta^1$ & 2.5 & 1.5 & 2.4512 & $1.95\times 10^{-2}$ \\
    $\zeta$ & 0.55 & 0.5 & 0.5508 & $1.45\times 10^{-3}$ \\
    \bottomrule
    \end{tabular}
    \caption{Test 4.1 (a): SIR transport model with spatially variable transmission rate and partially observed dynamics in diffusive regime ($\lambda_{S,I,R}^2=10^5$, $\tau_{S,I,R}=10^{-5}$). Inferred results for the three different coefficients in the incidence function, $\beta^0,\,\beta^1,\,\zeta$, and relative error with respect to the correct solution. }
    \label{Table:test4_parabolic_inference}
\end{table}
\begin{figure}[!t]
    \centering   
    \includegraphics[width=0.32\textwidth]{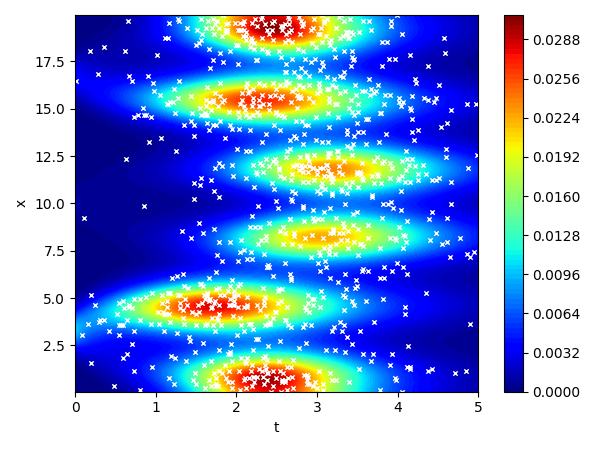}
    \includegraphics[width=0.32\textwidth]{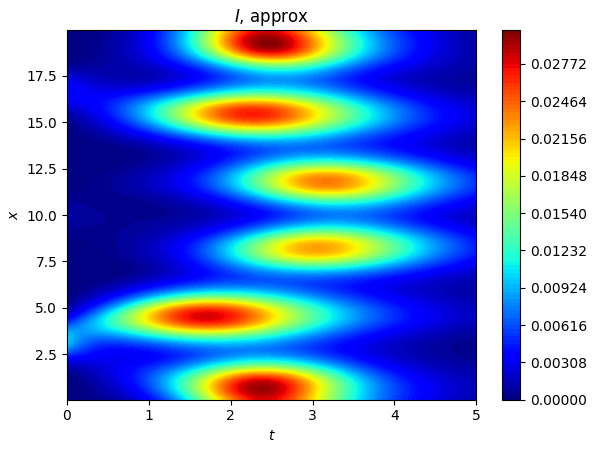}
    \includegraphics[width=0.32\textwidth]{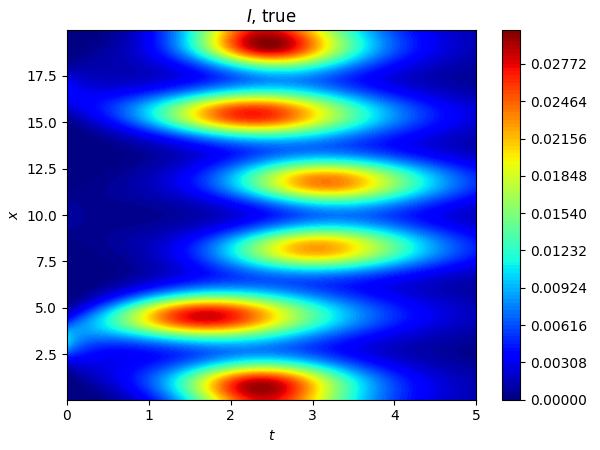}
    \caption{Test 4.1 (a): SIR transport model with spatially variable transmission rate and partially observed dynamics in diffusive regime ($\lambda_{S,I,R}^2=10^5$, $\tau_{S,I,R}=10^{-5}$). Selected sparse samples for the dataset marked with white crosses (left column), approximation obtained in the inverse problem (middle column), and ground truth (right column) of the densities of infected $I$. }
    \label{fig:test4_parabolic_inverse}
\end{figure}

\begin{figure}[!t]
    \centering   
    \includegraphics[width=0.32\textwidth]{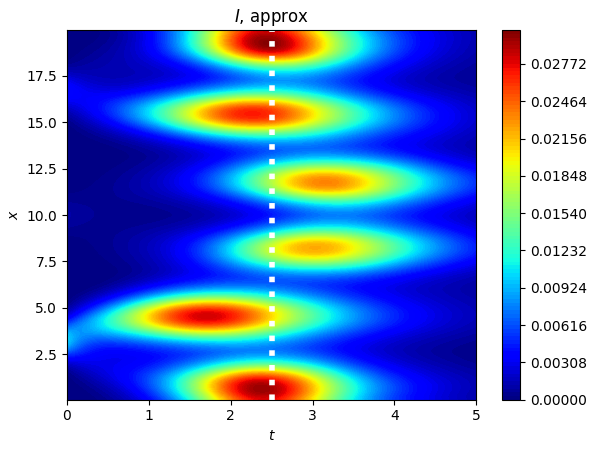}
    \includegraphics[width=0.32\textwidth]{figures_final/test3_new/scattered/test3_parabolic_scattered_I_true}
    \includegraphics[width=0.32\textwidth]{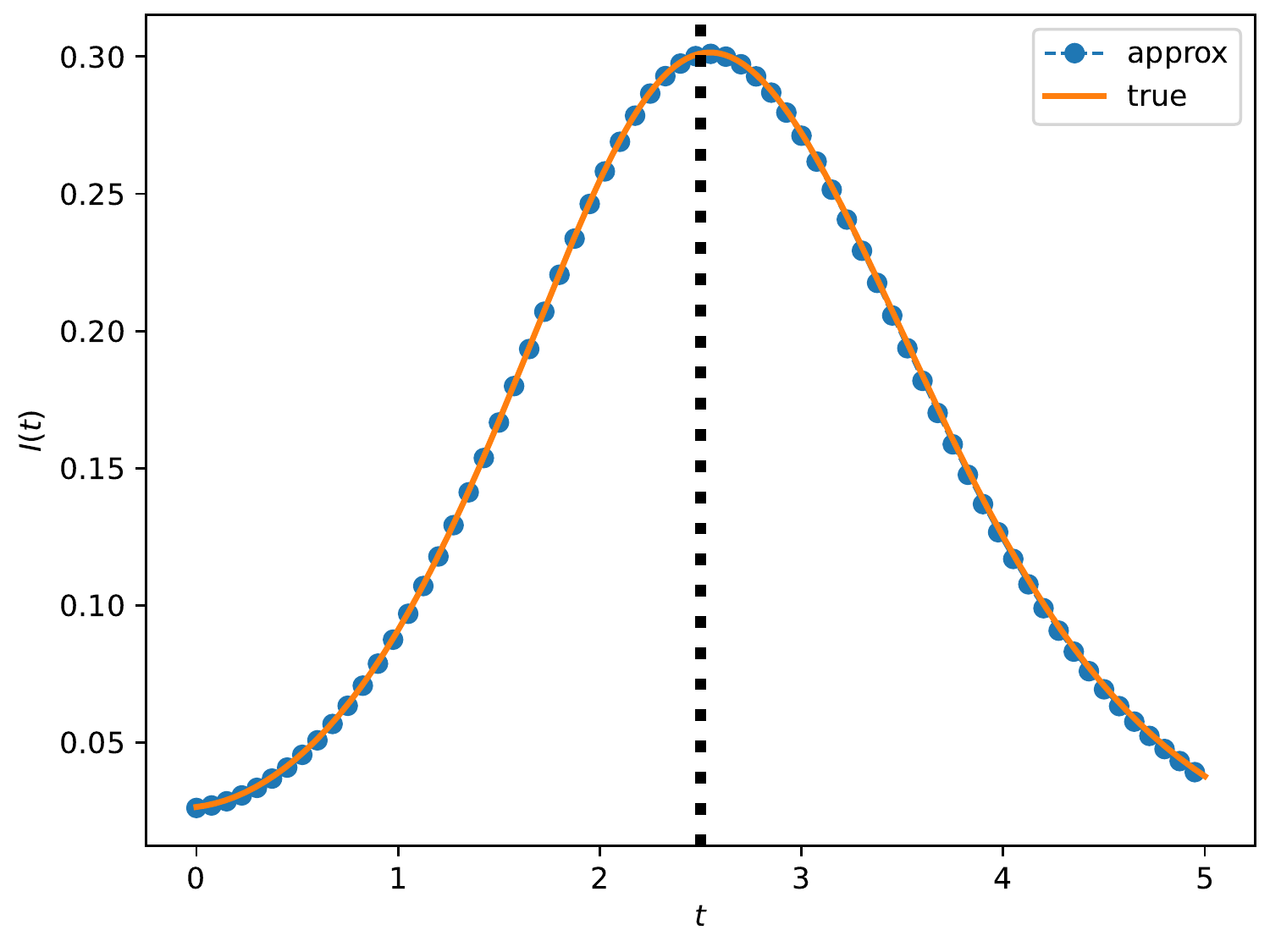}
    \includegraphics[width=0.32\textwidth]{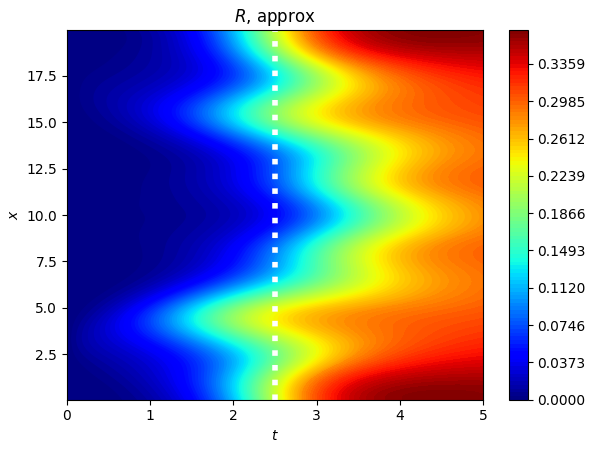}
    \includegraphics[width=0.32\textwidth]{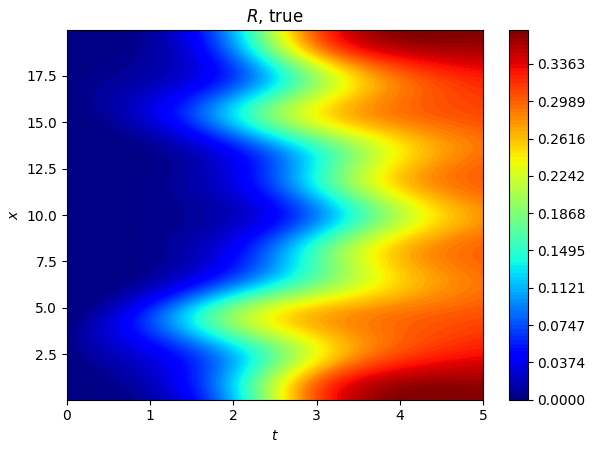}
    \includegraphics[width=0.32\textwidth]{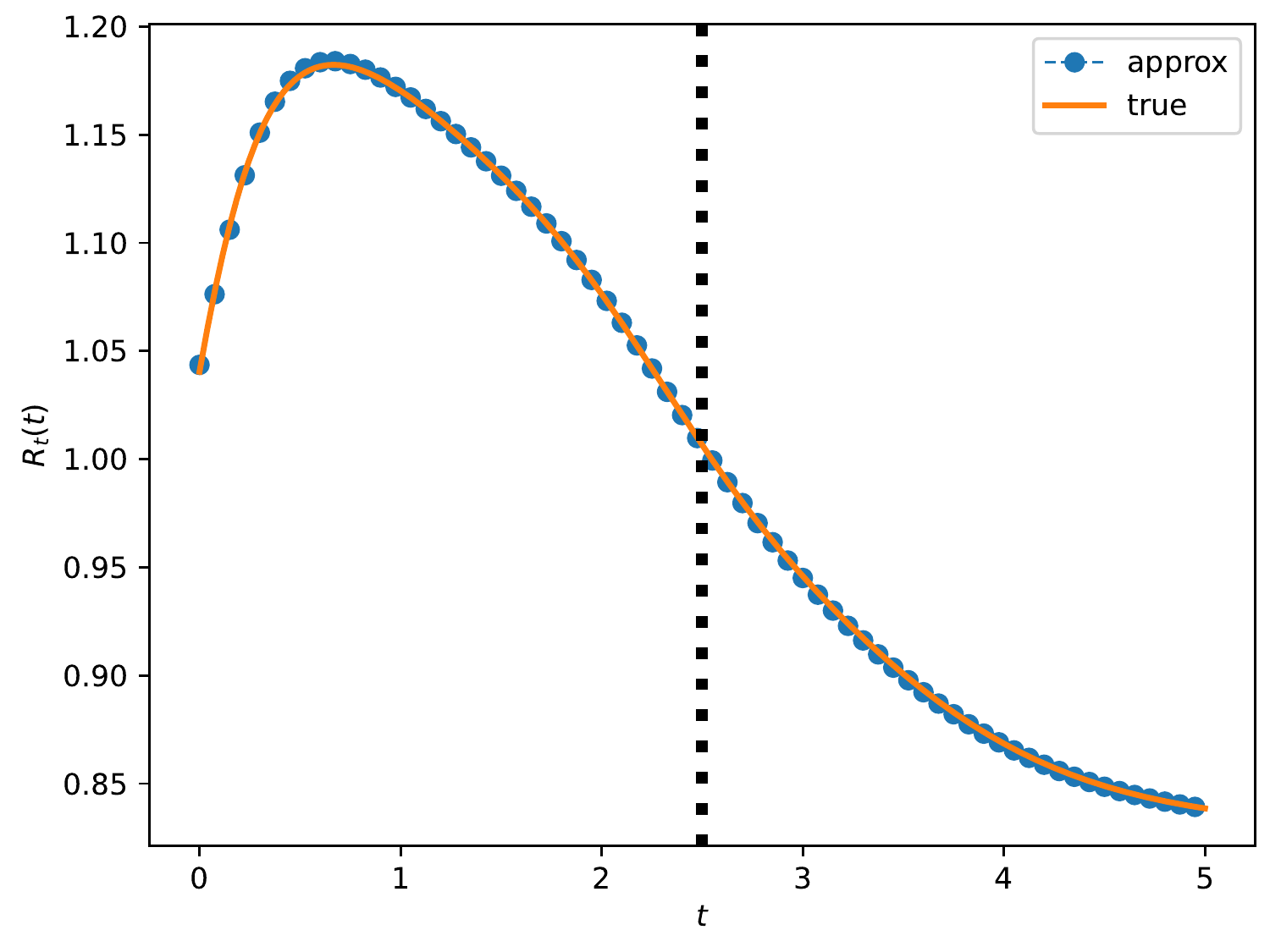}
    \caption{Test 4.2 (a): SIR transport model with spatially variable transmission rate and partially observed dynamics in diffusive regime ($\lambda_{S,I,R}^2=10^5$, $\tau_{S,I,R}=10^{-5}$). Approximation and forecast with measurements taken from a shorter time period, stopped at the dashed line (left column), and ground truth (middle column) of the densities of infected $I$ (first row) and removed $R$ (second row). Temporal evolution of the cumulative density of infected individuals $I$ in the whole domain (first row, right) and of the effective reproduction number $R_t$ (second row, right) obtained with the APNN, trained  based on measurements from a shorter time period $t \in [0,2.5]$ (marked by the dotted line). }
    \label{fig:test4_parabolic_forward}
\end{figure}
\subsubsection*{Test 4 (a): Partially observed dynamics in diffusive regime}
In the first scenario, a parabolic configuration of speeds and relaxation parameters is considered with $\lambda_{S,I,R}^2=10^5$ and $\tau_{S,I,R}=10^{-5}$. 

Similar to  the setting of Test 3,
we investigate the capabilities of the proposed APNN when concerning heterogeneous epidemic environments through the following two scenarios.
\begin{itemize}
\item {\bf Test 4.1 (a): Parameter inference test.} First, we consider a relatively sparse availability of measurements, with $N_d=1000$ samples over the spatio-temporal domain $(x, t)\in [0, 20]\times [0, 5]$  selected according to \eqref{eq:importance_sampling}, as previously described, with the main task to infer unknown physical parameters $\beta^0, \beta^1$ and $\zeta$. The selected samples are indicated in Figure \ref{fig:test4_parabolic_inverse} (left). 
\item {\bf Test 4.2 (a): Forecasting test.} Secondly, the forecasting performance in predicting the spread of the infectious disease until $t_{end} = 5$ with $N_d=10100$ measurements generated
from a shorter time period $t \in [0, 2.5]$ is investigated. 
\end{itemize}
In both scenarios, the equation residual is enforced on $N_r=10100$ residual points in the domain $(x, t)\in [0, 20]\times [0, 5]$, and initial conditions are enforced on $N_b=200$ equally spaced points. Furthermore, we enforce the conservation \eqref{eq:conservation} on $N_c=235$ equally spaced temporal points, and we randomly split 20\% of each dataset for validation purpose.

In Table \ref{Table:test4_parabolic_inference}, we present the results of parameters inference based on the sparse measurements. The APNN accurately recovers the correct values for parameters $\beta^0, \beta^1$ and $\zeta$ characterizing the epidemic incidence function, even when initial guesses are away from corresponding ground truth values. As illustrated in Figure \ref{fig:test4_parabolic_inverse}, reconstruction of the density $I$ is also in a very good agreement with the ground truth. Notice that the three initial epidemic concentrations give rise to six different epidemic outbreaks in time due to the spatial heterogeneity assigned to the transmission rate.
In Figure \ref{fig:test4_parabolic_forward}, we also present the approximated forward solutions for the forecasting task. A good match in the forecasting region ($t>2.5$) is observed, demonstrating once more the capability of the APNN to capture the underlying physics and deliver reasonably accurate predictions in the forecasting regions, even when spatially heterogeneous environments are considered in the context of partially observed systems. 

\begin{table}[!b]
    \centering
    \begin{tabular}{@{}ccccc@{}}\toprule 
    Parameter & Ground Truth & Initial Guess & Estimation & Relative Error \\
    \midrule 
    $\beta^0$ & 9 & 5 & 9.0205 & $2.28\times 10^{-3}$ \\
    $\beta^1$ & 2.5 & 1.5 & 2.4691 & $1.24\times 10^{-2}$ \\
    $\zeta$ & 0.55 & 0.5 & 0.5502 & $3.64\times 10^{-4}$ \\
    \bottomrule
    \end{tabular}
    \caption{Test 4.1 (b): SIR transport model with spatially variable transmission rate and partially observed dynamics in hyperbolic regime ($\lambda_{S,I,R}^2=1$, $\tau_{S,I,R}=1$). Inferred results from sparse measurements for the three different coefficients in the incidence function, $\beta^0,\,\beta^1,\,\zeta$, and relative error with respect to the ground truth.}
    \label{Table:test4_hyperbolic_inference}
\end{table}
\begin{figure}[!t]
    \centering
    \includegraphics[width=0.32\textwidth]{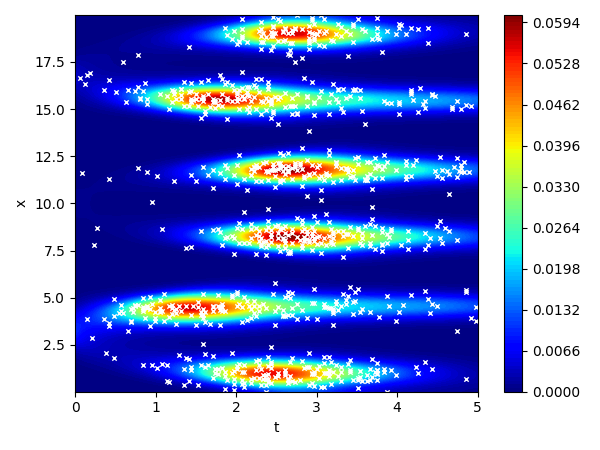}   
    \includegraphics[width=0.32\textwidth]{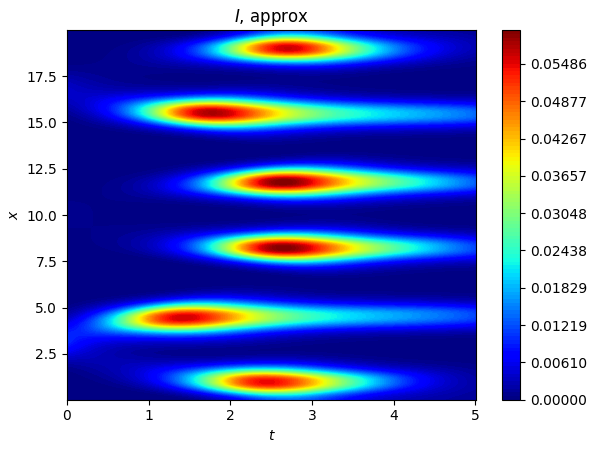}
    \includegraphics[width=0.32\textwidth]{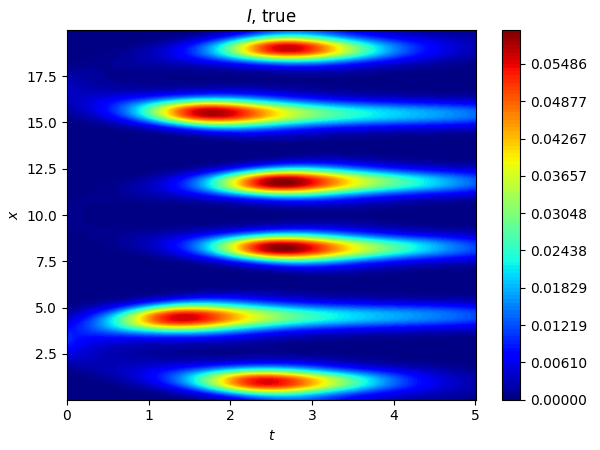}
    \caption{Test 4.1 (b): SIR transport model with spatially variable transmission rate and partially observed dynamics in hyperbolic regime ($\lambda_{S,I,R}^2=1$, $\tau_{S,I,R}=1$). Selected sparse samples for the dataset marked with white crosses (left column), approximation obtained in the inverse problem (middle column), and ground truth (right column) of the densities of infected $I$.}
    \label{fig:test4_hyperbolic_inverse}
\end{figure}

\begin{figure}[!t]
    \centering   
    \includegraphics[width=0.32\textwidth]{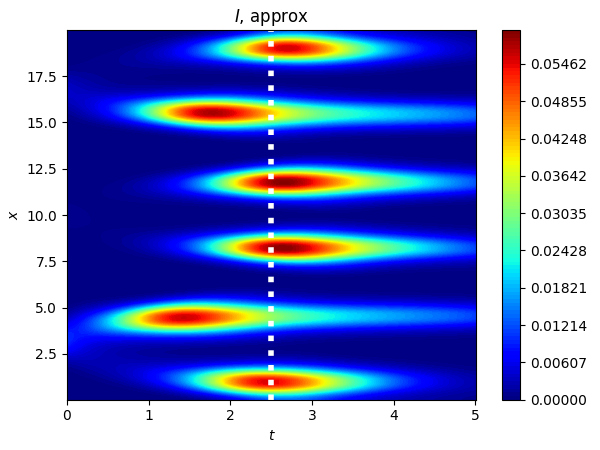}
    \includegraphics[width=0.32\textwidth]{figures_final/test3_new/scattered/test3_hyperbolic_scattered_I_true}
    \includegraphics[width=0.32\textwidth]{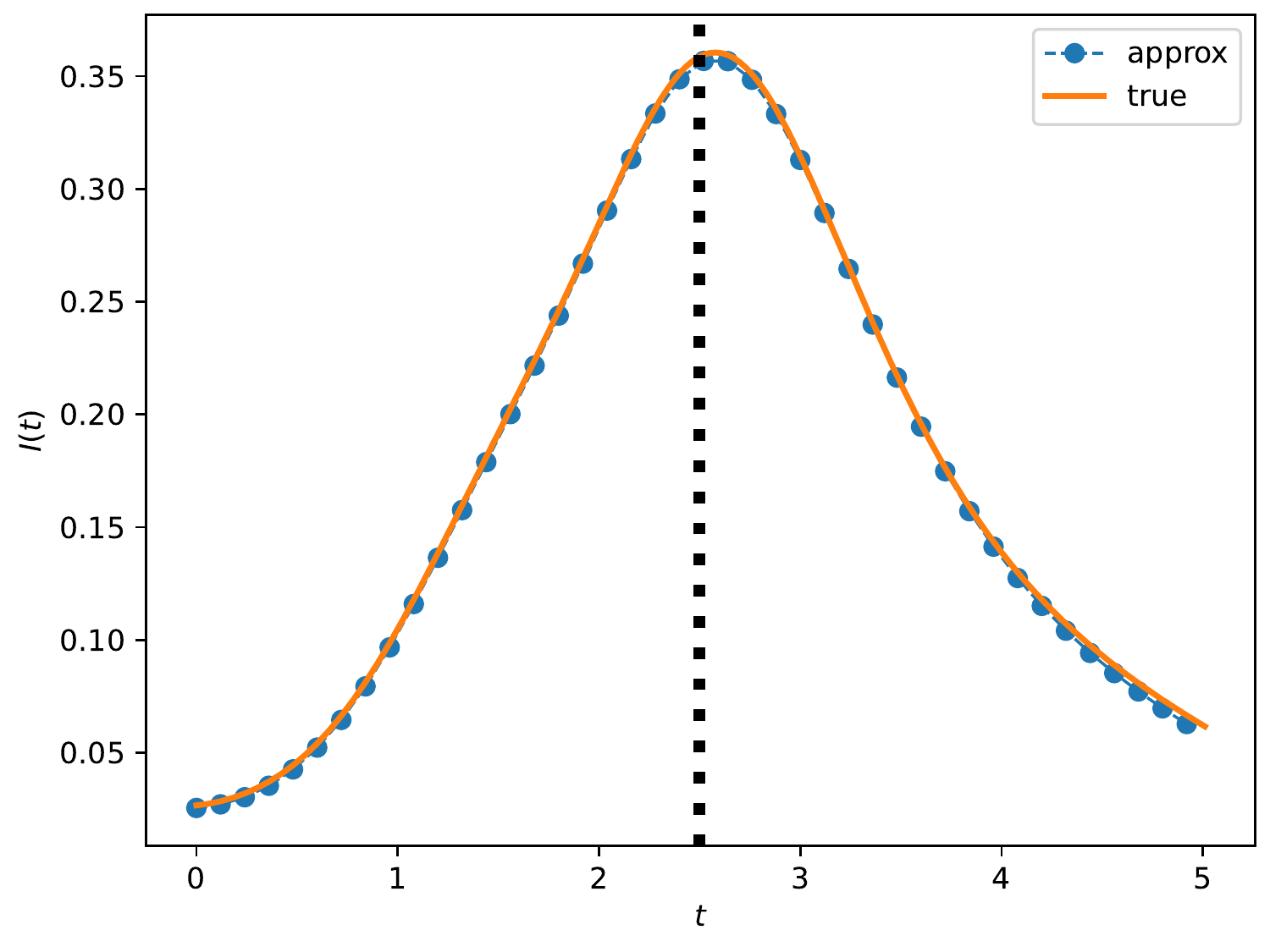}
    \includegraphics[width=0.32\textwidth]{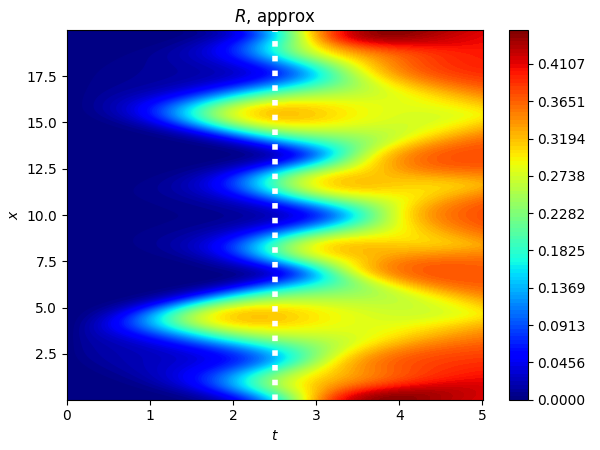}
    \includegraphics[width=0.32\textwidth]{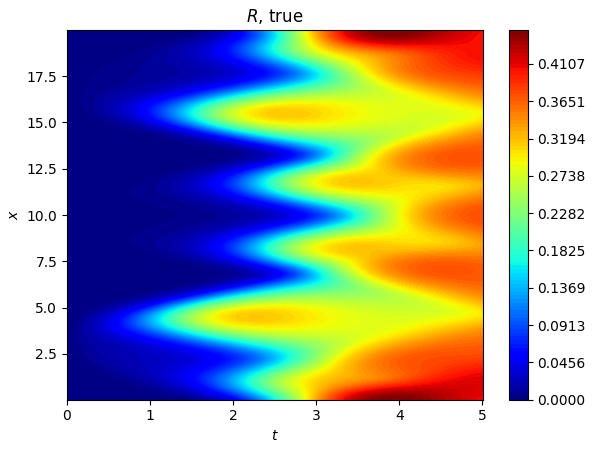}
    \includegraphics[width=0.32\textwidth]{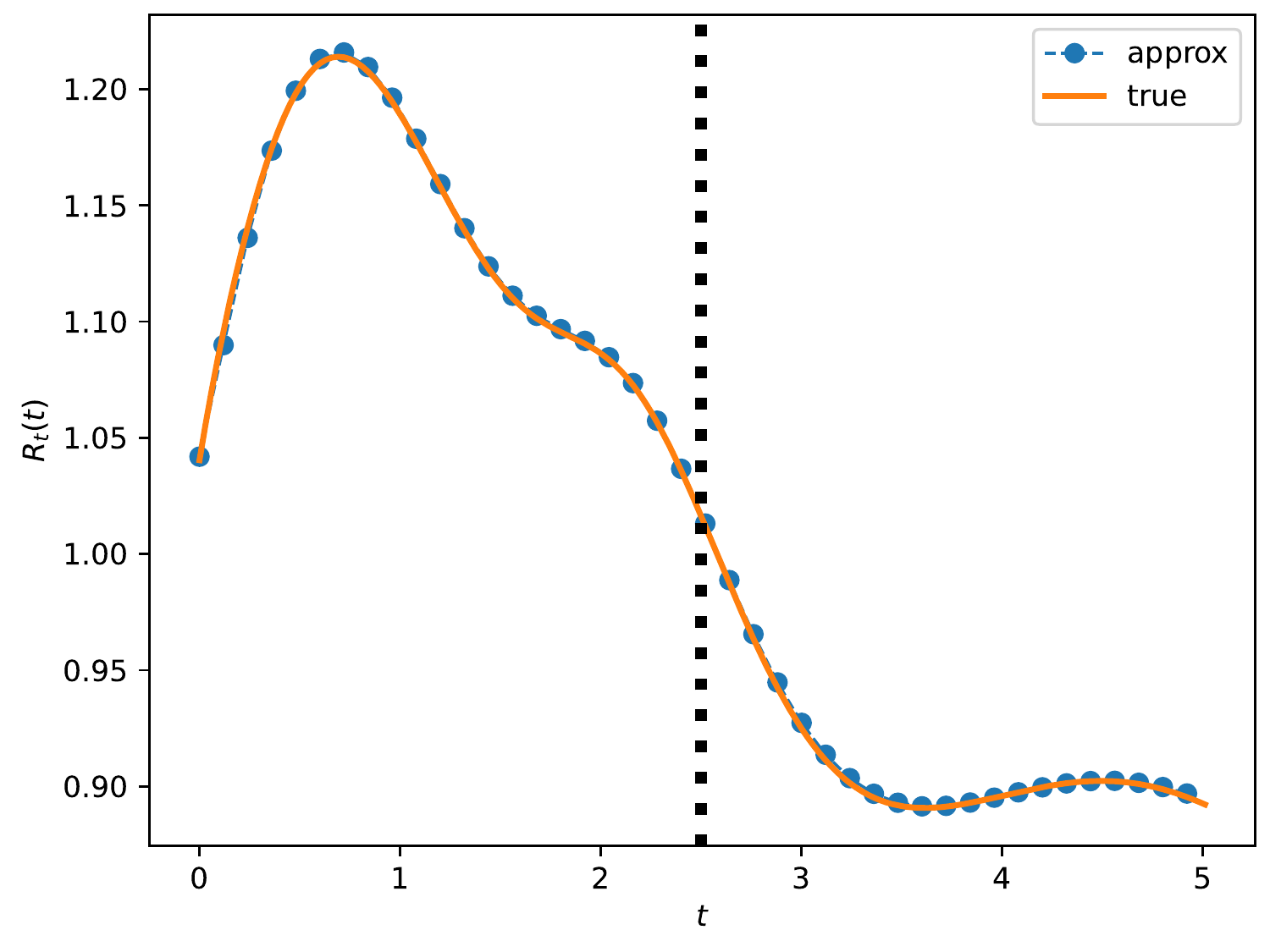}
    \caption{Test 4.2 (b): SIR transport model with spatially variable transmission rate and partially observed dynamics in hyperbolic regime ($\lambda_{S,I,R}^2=1$, $\tau_{S,I,R}=1$). Approximation and forecast with measurements taken from a shorter time period, stopped at the dashed line (left column), and ground truth (middle column) of the densities of infected $I$ (first row) and removed $R$ (second row). Temporal evolution of the cumulative density of infected individuals $I$ in the whole domain (first row, right) and of the effective reproduction number $R_t$ (second row, right) obtained with the APNN, trained  based on measurements from a shorter time period $t \in [0,2.5]$ (marked by the dotted line). }
    \label{fig:test4_hyperbolic_forward}
\end{figure}

\subsubsection*{Test 4 (b): Partially observed dynamics in hyperbolic regime}
In the second scenario, we consider a hyperbolic regime setting $\lambda_{S,I,R}=1$ and $\tau_{S,I,R}=1$. Similar to the previous test case, we consider two distinguished tasks for the APNN. 
\begin{itemize}
\item {\bf Test 4.1 (b): Parameter inference test.} Initially, a sparse measurement dataset of $N_d=1000$ training samples over the spatio-temporal domain $(x, t)\in [0, 20]\times [0, 5]$ is considered, based on the importance sampling previously described, and marked in Figure \ref{fig:test4_hyperbolic_inverse} (left), to solve the inverse problem and also evaluate the following forward reconstruction. 
\item {\bf Test 4.2 (b): Forecasting test.} Then, the APNN is trained with $N_d=8400$ data samples selected from the spatio-temporal domain $(x, t)\in [0, 20]\times [0, 2.5]$, and the reconstruction of the dynamics is evaluated until $t=5$, to also examine the performance on the forecasting of the virus spread.
\end{itemize}
Equation residual is enforced on $N_r=23500$ residual points on the domain $(x, t)\in [0, 20]\times [0, 5]$ in both setups, while $N_b=600$ points are applied to enforce initial conditions for densities $S,I,R$ and fluxes $J_S,J_I,J_R$, and the conservation \eqref{eq:conservation} is enforced on $N_c=47$ equally spaced temporal points. 20\% of each dataset is randomly selected, as usual, for validation during the training process. 

Similarly with the parabolic setting, the APNN is able to estimate the correct parameters of the spatially-varied transmission rate $\beta$ from sparse measurements in the hyperbolic regime, as shown in Table \ref{Table:test4_hyperbolic_inference}. In Figures \ref{fig:test4_hyperbolic_inverse} and \ref{fig:test4_hyperbolic_forward}, the approximated forward solutions and the ground truth of $I$ and $R$ over the space-time domain are shown. A good match between the APNN approximation and the ground truth is observed, for both sparse measurement and the measurement from a reduced training time domain, in the latter considering also predictions of the space-time dynamics. Notice that, as expected, due to the hyperbolic setting of the scaling parameters of this test, the six epidemic outbreaks that arise at different temporal levels due to the spatial movement of individuals are more contained in terms of spatial spread with respect to results obtained in the diffusive regime.

\section{Conclusions}
The recent Covid-19 pandemic has led to a significant development of mathematical models for describing epidemiological phenomena, which have also introduced the challenge of identifying the parameters involved from partial information. In this direction, recent developments in machine learning represent a promising tool for addressing such problems in the hope of identifying robust procedures for solving the corresponding inverse problems and also formulating predictive scenarios. This paper has addressed these problems in the context of spatially dependent epidemic models for which, in addition to the lack of information about the spread of the epidemic, face additional difficulties induced by the different scales at which the dynamics take place. These scales are representative of the different interactions that occur in densely populated areas, such as urban areas, or in suburban areas where the movement of individuals over long distances prevails. The construction of neural networks that can accurately describe the various scales is thus essential. In particular, we have shown how physically informed neural networks (PINN) that benefit from the asymptotic-preserving (AP) property provide considerably better results with respect to the different scales of the problem when compared with standard PINN.   
Several numerical tests have been presented to illustrate the performance of this new class of neural networks, referred to as asymptotic-preserving neural network (APNN), both for inverse and forward problems. Finally, we emphasize that even if, for presentation simplicity, we focused on a single population hyperbolic SIR model, the results extend naturally to multi-population transport models which include additional epidemic compartments^^>\cite{ABBDPTZ21, Bert3, BDP}.

\section*{Acknowledgments}
G.B. and  L.P. were partially supported by MIUR-PRIN Project 2017, No. 2017KKJP4X \emph{Innovative numerical methods for evolutionary partial differential equations and applications}. G.B. also acknowledges the support by INdAM--GNCS. X.Z. was supported by the Simons Foundation (504054). 
\appendix

\section{AP-loss function for the Goldstein--Taylor model}
\label{appendix}
Fixing a finite set of residual points $\{(x_r^n, t_r^n)\}_{n=1}^{N_r} \subset \Omega$, $\{(x_b^k, t_b^k)\}_{k=1}^{N_b} \subset \partial \Omega$, and considering the available dataset $\{u_d^i, x_d^i, t_d^i \}_{i=1}^{N_d}$, we define the loss function for the Goldstein--Taylor model as follows 
\begin{equation}
    \label{eq:pinn-loss_GT}
    \begin{split}
    \mathcal{L}(\theta) = &\, \omega_d^T \,{\mathcal L}_d(\theta) + \omega_b^T \,{\mathcal L}_b(\theta) + \omega_r^T \,{\mathcal L}_r^{\epsilon}(\theta).
    \end{split}
\end{equation}

The expressions of ${\mathcal L}_d$ and ${\mathcal L}_b$ terms in \eqref{eq:pinn-loss_GT} read
\begin{equation}
\begin{split}
    \omega_d^T\, {\mathcal L}_d(\theta)  
    = &\frac{\omega_d^\rho}{N_d} \sum_{i=1}^{N_d} \left| \rho_{NN}(x_d^i, t_d^i;\theta) - \rho(x_d^i, t_d^i) \right|^2+
    \frac{\omega_d^j}{N_d} \sum_{i=1}^{N_d} \left| j_{NN}(x_d^i, t_d^i;\theta) - j(x_d^i, t_d^i) \right|^2\,,\\
    \omega_b^T\, {\mathcal L}_b(\theta)  
    = &\frac{\omega_b^\rho}{N_b} \sum_{k=1}^{N_b} \left| \rho_{NN}(x_b^k, t_b^k;\theta) - \rho(x_b^k, t_b^k) \right|^2+
    \frac{\omega_b^j}{N_b} \sum_{k=1}^{N_b} \left| j_{NN}(x_b^k, t_b^k;\theta) - j(x_b^k, t_b^k) \right|^2\,.
\end{split}
\label{eq:pinn-MSEd_GT}
\end{equation}

The AP-residual term ${\mathcal L}_r^{\epsilon}$ in \eqref{eq:pinn-loss_GT} is defined through \eqref{eq:pinn-ap-residue_GT} and reads as follows
\begin{equation}
\begin{split}
    \omega_r^T\, {\mathcal L}_r^{\epsilon}(\theta)
    = &\frac{\omega_r^\rho}{N_r} \sum_{n=1}^{N_r} \left| \frac{\partial \rho_{NN}(x_r^n, t_r^n;\theta)}{\partial t} + \frac{\partial j_{NN}(x_r^n, t_r^n;\theta)}{\partial x} \right|^2\\+
    &\frac{\omega_r^j}{N_r} \sum_{n=1}^{N_r} \left| \epsilon^2 \frac{\partial j_{NN}(x_r^n, t_r^n;\theta)}{\partial t} + c^2\frac{\partial \rho_{NN}(x_r^n, t_r^n;\theta)}{\partial x} + \sigma\, j_{NN}(x_r^n, t_r^n;\theta) \right|^2,
\end{split}
\label{eq:pinn-MSEr_GT}
\end{equation}
where we assumed that the scattering coefficient $\sigma$ in system \eqref{Goldstein-Taylor} is constant.

This formulation allows the neural network to benefit from the AP property. Indeed, if we consider the zero relaxation (or diffusive) limit $\epsilon\to 0$, Eq. \eqref{eq:pinn-MSEr_GT} results
\begin{equation*}
\begin{split}
    \omega_r^T{\mathcal L}_r^{0}(\theta)
    = &\frac{\omega_r^\rho}{N_r} \sum_{n=1}^{N_r} \left| \frac{\partial \rho_{NN}(x_r^n, t_r^n;\theta)}{\partial t} + \frac{\partial j_{NN}(x_r^n, t_r^n;\theta)}{\partial x} \right|^2 \\+
    &\frac{\omega_r^j}{N_r} \sum_{n=1}^{N_r} \left| c^2\frac{\partial \rho_{NN}(x_r^n, t_r^n;\theta)}{\partial x} + \sigma_{NN}\, j_{NN}(x_r^n, t_r^n;\theta) \right|^2\,,
\end{split}
\label{eq:pinn-MSEr_GT_diff}
\end{equation*}
which is consistent with the definition of residual term that we would have by directly considering the diffusive limit of the model in \eqref{eq.heat}.

\section{AP-loss function for the hyperbolic SIR model}
\label{appendix1}
Aware that in real-world applications data on fluxes $J_S,J_I,J_R$ are generally difficult to access, we consider the loss function 
\begin{equation}
\mathcal{L}(\theta) = \omega_d^T \,{\mathcal L}_d(\theta) + \omega_b^T\, {\mathcal L}_b(\theta) + \omega_r^T \,{\mathcal L}^\tau_r(\theta) + \omega_c\, {\mathcal L}_c(\theta), 
\label{eq:pinn-loss2}
\end{equation}
where the first term is given by
\begin{equation*}
\begin{split}
    \omega_d^T \,{\mathcal L}_d(\theta) =
    &\frac{\omega_d^S}{N_d} \sum_{i=1}^{N_d} \left| S_{NN}(x_d^i, t_d^i;\theta) - S(x_d^i, t_d^i) \right|^2 +
    \frac{\omega_d^I}{N_d} \sum_{i=1}^{N_d} \left| I_{NN}(x_d^i, t_d^i;\theta) - I(x_d^i, t_d^i) \right|^2 \\+
    &\frac{\omega_d^R}{N_d} \sum_{i=1}^{N_d} \left| R_{NN}(x_d^i, t_d^i;\theta) - R(x_d^i, t_d^i) \right|^2\,.
\end{split}
\end{equation*} 
Concerning the imposition of initial and boundary conditions, both included in the ${\mathcal L}_b$ term in \eqref{eq:pinn-loss2}, we have the detailed expression 
\begin{equation*}
\begin{split}
    \omega_b^T \,{\mathcal L}_b(\theta) =
    &\frac{\omega_b^S}{N_b} \sum_{k=1}^{N_b} \left| S_{NN}(x_b^k, t_b^k;\theta) - S(x_b^k, t_b^k) \right|^2 +
    \frac{\omega_b^I}{N_b} \sum_{k=1}^{N_b} \left| I_{NN}(x_b^k, t_b^k;\theta) - I(x_b^k, t_b^k) \right|^2 \\+
    &\frac{\omega_b^R}{N_b} \sum_{k=1}^{N_b} \left| R_{NN}(x_b^k, t_b^k;\theta) - R(x_b^k, t_b^k) \right|^2 +
    \frac{\omega_b^{J_S}}{N_b} \sum_{k=1}^{N_b} \left| J^S_{NN}(x_b^k, t_b^k;\theta) - J_S(x_b^k, t_b^k) \right|^2 \\+
    &\frac{\omega_b^{J_I}}{N_b} \sum_{k=1}^{N_b} \left| J^I_{NN}(x_b^k, t_b^k;\theta) - J_I(x_b^k, t_b^k) \right|^2 +
    \frac{\omega_b^{J_R}}{N_b} \sum_{k=1}^{N_b} \left| J^R_{NN}(x_b^k, t_b^k;\theta) - J_R(x_b^k, t_b^k) \right|^2\,.
\end{split}
\end{equation*} 
We underline here that, to impose boundary conditions, we apply the appropriate mapping technique based on the specific boundary conditions of the problem of interest^^>\cite{zhang2020learning}. 

The expression of the last mean squared error term in the loss function \eqref{eq:pinn-loss2}, which concerns the residual presented in \eqref{eq:pinn-ap-residue}, reads
\begin{alignat*}{2}
\omega_r^T\, {\mathcal L}_r^\tau(\theta) = &\frac{\omega_r^S}{N_r} \sum_{n=1}^{N_r} &&\bigg|\frac{\partial S_{NN}(x_r^n, t_r^n;\theta)}{\partial t} + \frac{\partial J^S_{NN}(x_r^n, t_r^n;\theta)}{\partial x}
+ \beta(x_r^n, t_r^n) S_{NN}(x_r^n, t_r^n;\theta) I_{NN}(x_r^n, t_r^n;\theta)\bigg|^2 \\
+ &\frac{\omega_r^I}{N_r} \sum_{n=1}^{N_r} &&\bigg| \frac{\partial I_{NN}(x_r^n, t_r^n;\theta)}{\partial t} + \frac{\partial J^I_{NN}(x_r^n, t_r^n;\theta)}{\partial x}
- \beta(x_r^n, t_r^n) S_{NN}(x_r^n, t_r^n;\theta) I_{NN}(x_r^n, t_r^n;\theta) \\
& &&+ \gamma(x_r^n, t_r^n) I_{NN}(x_r^n, t_r^n;\theta) \bigg|^2 \\
+ &\frac{\omega_r^R}{N_r} \sum_{n=1}^{N_r} &&\bigg| \frac{\partial R_{NN}(x_r^n, t_r^n;\theta)}{\partial t} + \frac{\partial J^R_{NN}(x_r^n, t_r^n;\theta)}{\partial x} -  \gamma(x_r^n, t_r^n) I_{NN}(x_r^n, t_r^n;\theta)  \bigg|^2 \\+
&\frac{\omega_r^{J_S}}{N_r} \sum_{n=1}^{N_r} &&\bigg| \tau_S(x_r^n) \frac{\partial J^S_{NN}(x_r^n, t_r^n;\theta)}{\partial t} + D_S(x_r^n) \frac{\partial S_{NN}(x_r^n, t_r^n;\theta)}{\partial x} \\
& &&+ \tau_S(x_r^n) \beta(x_r^n, t_r^n) J^S_{NN}(x_r^n, t_r^n;\theta) I_{NN}(x_r^n, t_r^n;\theta) + J^S_{NN}(x_r^n, t_r^n;\theta)  \bigg|^2\\+
&\frac{\omega_r^{J_I}}{N_r} \sum_{n=1}^{N_r} &&\bigg| \tau_I(x_r^n) \frac{\partial J^I_{NN}(x_r^n, t_r^n;\theta)}{\partial t} + D_I(x_r^n) \frac{\partial I_{NN}(x_r^n, t_r^n;\theta)}{\partial x} \\
& &&- \tau_I(x_r^n) \frac{\lambda_I(x_r^n)}{\lambda_S(x_r^n)} \beta(x_r^n, t_r^n) J^S_{NN}(x_r^n, t_r^n;\theta) I_{NN}(x_r^n, t_r^n;\theta) \\
& &&+ \tau_I(x_r^n) \gamma(x_r^n, t_r^n) J^I_{NN}(x_r^n, t_r^n;\theta) + J^I_{NN}(x_r^n, t_r^n;\theta)  \bigg|^2\\+
&\frac{\omega_r^{J_R}}{N_r} \sum_{n=1}^{N_r} &&\bigg| \tau_R(x_r^n) \frac{\partial J^R_{NN}(x_r^n, t_r^n;\theta)}{\partial t} + D_R(x_r^n) \frac{\partial R_{NN}(x_r^n, t_r^n;\theta)}{\partial x} \\ 
& &&- \tau_R(x_r^n) \frac{\lambda_R(x_r^n)}{\lambda_I(x_r^n)} \gamma(x_r^n, t_r^n) J^I_{NN}(x_r^n, t_r^n;\theta) + J^R_{NN}(x_r^n, t_r^n;\theta) \bigg|^2\,.
\end{alignat*}
If one considers the limit $\tau_i\to 0$, $\lambda_i\to \infty$, $i\in\{S,I,R\}$ s.t. equations \eqref{eq:diff2} hold true, the above term is clearly consistent with the definition of residual applied directly to the diffusive limit given by equations \eqref{eq:diff}.

Finally, in the SIR epidemic transport model, conservation \eqref{eq:conservation} is enforced as
\begin{equation}
    \begin{aligned}
        \omega_c\, {\mathcal L}_c(\theta) = \frac{\omega_c}{N_c}\sum_{m=1}^{N_c} &\left|\sum_{q=1}^{N_q} \left(S_{NN}(x_q, t_c^m;\theta) + I_{NN}(x_q, t_c^m;\theta) + R_{NN}(x_q, t_c^m;\theta) \right)\right. \\
        & \left.-\sum_{q=1}^{N_q} \left(S(x_q, 0) + I(x_q, 0) + R(x_q, 0)\right)\right|^2,
    \end{aligned}
\end{equation}
with $N_q$ quadrature points in $\mathcal{D}$.


\section{Loss function and training hyperparameters}
\label{appendix_C}

Loss function and training hyperparameters of the APNN for the various test cases considered are listed in the following Tables.

\label{sec:appendix-model-hyperparameters}
\begin{table}[!h]
    \centering
    \begin{tabular}{@{}lccccc@{}}\toprule 
    Test &$\omega_d^\rho$ & $(\omega_b^{\rho}, \omega_b^j)$ & $(\omega_r^\rho, \omega_r^j)$ & LR & Epochs \\
    \midrule 
    1 &100 & (1, 1) & (1, 1) & $10^{-2}$ & 20000 \\
    2 (a) &100 & (1, 1) & (1, 1) & $10^{-2}$ & 40000 \\
    2 (b) &100 & (1, 1) & (1, 1) & $10^{-2}$ & 12000 \\
    \bottomrule
    \end{tabular}
    \caption{Goldstein-Taylor model. Loss function and training hyperparameters of the APNN. LR refers to the learning rate used in the optimization. In all the tests we consider full batch size.}
    \label{Table:test12_hyperparameters}
\end{table}

\begin{table}[!h]
{\small
    \centering
    \begin{tabular}{@{}lccccc@{}}\toprule 
    Test &$(\omega_u^S,\,\omega_u^I,\,\omega_u^R)$ &$(\omega_b^S,\,\omega_b^I,\,\omega_b^R,\,\omega_b^{J_S},\,\omega_b^{J_I},\,\omega_b^{J_R})$ &$(\omega_r^S,\,\omega_r^I,\,\omega_r^R,\,\omega_r^{J_S},\,\omega_r^{J_I},\,\omega_r^{J_R})$ &LR & Epochs \\
    \midrule 
    3 (a) &(1, 100, 10) & (1, 10, 1,-,-,-) & (1, 10, 1, 1, 10, 1) & $10^{-3}$ & 20000\\
    3 (b) &(1, 100, 10) & (1, 10, 1,-,-,-) & (1, 100, 10, 1, 100, 10) & $10^{-2}$ & 20000\\
    4 (a) &(1, 1000, 100) & (1, 1000, 100, 1, 10, 1) & (1, 1000, 100, 1, 10, 1) & $10^{-3}$ & 150000\\
    4 (b) &(1, 1000, 100) & (1, 1000, 100, 1, 10, 1) & (1, 1000, 100, 1, 10, 1) & $10^{-3}$ & 150000\\
    \bottomrule
    \end{tabular}}
    \caption{SIR hyperbolic model. Loss function and training hyperparameters of the APNN. We fix the weight for the conservation loss in SIR models $\omega_c=1$. LR is the learning rate used in the optimization. In all the tests we consider full batch size. }
    \label{Table:test34_hyperparameters}
\end{table}


\end{document}